\documentclass[a4 paper,12pt]{amsart}
\usepackage{amscd}
\usepackage{amssymb}

\usepackage{graphicx}
\usepackage{comment}
\usepackage[margin=2 in]{geometry}
\usepackage{setspace}
\usepackage{enumitem}

\usepackage[all]{xy}

\usepackage{color}
\allowdisplaybreaks[1]

\newtheorem{thm}{Theorem}[section]

\newtheorem{cor}[thm]{Corollary}
\newtheorem{lem}[thm]{Lemma}
\newtheorem{prop}[thm]{Proposition}

\newtheorem{notation}[thm]{Notation}

\theoremstyle{definition}

\numberwithin{equation}{section}

\setlist[itemize]{leftmargin=*}

\def\O{\Omega}

\newcommand{\C}{\mathbb C}

\newcommand{\R}{\mathbb R}

\newcommand{\N}{\mathbb N}

\newcommand{\HD}[1]{\mathrm{H}^{#1}(\mathfrak{sl}(2),\mathcal D_{{\lambda},\mu})}
\newcommand{\ZD}[1]{Z^{#1}(\mathfrak{sl}(2),\mathcal D_{{\lambda},\mu})}
\newcommand{\BD}[1]{B^{#1}(\mathfrak{sl}(2),\mathcal D_{{\lambda},\mu})}

\newcommand{\HDif}[1]{\mathrm{H}^{#1}_{\mathrm{diff}}(\mathfrak{sl}(2),\mathcal D_{{\lambda},\mu})}
\newcommand{\ZDif}[1]{Z^{#1}_{\mathrm{diff}}(\mathfrak{sl}(2),\mathcal D_{{\lambda},\mu})}
\newcommand{\BDif}[1]{B^{#1}_{\mathrm{diff}}(\mathfrak{sl}(2),\mathcal D_{{\lambda},\mu})}

\newcommand{\HM}[1]{\mathrm{H}^{#1}(\mathfrak{sl}(2),M)}
\newcommand{\ZM}[1]{Z^{#1}(\mathfrak{sl}(2),M)}
\newcommand{\BM}[1]{B^{#1}(\mathfrak{sl}(2),M)}

\newcommand{\Zred}[1]{Z^{#1}_{\mathrm{red}}(\mathfrak{sl}(2),M)}
\newcommand{\Bred}[1]{B^{#1}_{\mathrm{red}}(\mathfrak{sl}(2),M)}

\newcommand{\Dlm}{\mathcal D_{{\lambda},\mu}}

\newcommand{\coker}{\mathrm{coker}}

\def \l{\lambda}
\def \a{\alpha}

\def\1{\mathbf 1}

\newcommand{\Span}{{\rm Span}}

\newcommand{\resumename}{R\'esum\'e}

\usepackage{diagbox}

\begin{document} 

%----------------------------------------------------------------------------------------------------------------------

\title[$\mathfrak{sl}(2)$-cohomology]{ $\mathfrak{sl}(2)$-cohomology on multidifferential operators space}

\author[Arnal]{Didier Arnal}

\address{
Institut de Mathématiques de Bourgogne\\
UMR CNRS 5584\\
Université de Bourgogne Franche Comté\\
U.F.R. Sciences et Techniques
B.P. 47870\\
F-21078 Dijon Cedex\\France} \email{Didier.Arnal@wanadoo.fr}

\author[Ben Ammar]{Mabrouk Ben Ammar}

\address{
Laboratoire Alg\`ebre, G\'eom\'etrie et Th\'eorie Spectrale (AGTS), LR11ES53\\
Facult\'e des Sciences de Sfax,\\
Universit\'e de Sfax,\\
Route de Soukra, km 3,5, B.P. 1171, 3000 Sfax, Tunisie} \email{mabrouk.benammar@gmail.com}



%-----------------------------------------------------------------------------------------------------------------------

\begin{abstract}
In this paper, we consider the known question of the computation of the cohomology spaces for the natural action of the Lie algebra $\mathfrak{sl}(2)$ acting on the space $\Dlm$ of multidifferential operators on the line. We first show that these cohomology spaces are entirely characterized by a linear map $f$ between two explicit finite dimensional spaces, then we prove that this map has always maximal rank, allowing to get explicit recursive expression for their dimensions.\\
\end{abstract}

%-----------------------------------------------------------------------------------------------------------------------

\maketitle

%-----------------------------------------------------------------------------------------------------------------------

\section{Introduction}\label{introduction}

%-----------------------------------------------------------------------------------------------------------------------

\

In this paper, we consider the Lie algebra of vector fields on the real line:
$$
\mathrm{Vect}(\mathbb{R})=\Big\{X_h=h\frac{d}{dx}: ~ h\in C^\infty(\mathbb{R})\Big\},
$$
with its usual bracket $[X_{h_1},X_{h_2}]=X_{h_1h'_2-h_2h'_1}$. This Lie algebra acts naturally on the space of $\l$-densities on the line ($\l\in\R$), denoted by:
\begin{equation*}
\mathcal{F}_\l=\big\{ fdx^{\l}: ~f\in C^\infty(\mathbb{R})\big\}.
\end{equation*}
The action of $\mathrm{Vect}(\mathbb{R})$ on this space $\mathcal{F}_\l$ by Lie derivative is:
\begin{equation}\label{Lie1}
X_h\cdot (fdx^{\l}) =L_{X_h}^\l(fdx^{\l}):=(hf'+\l h'f)dx^\l.
\end{equation}

More generally, if $\l=(\l_1,\ldots,\l_n)$ is in $\R^n$, and $\mu$ in $\R$, we denote by $\Dlm$ the space of $n$-ary differential operators $A$ from $\mathcal{F}_{\l_1}\otimes\cdots\otimes\mathcal{F}_{\l_n}$ to $\mathcal{F}_\mu$. The Lie algebra $\mathrm{Vect}(\mathbb{R})$ acts on the space $\mathcal{D}_{{\lambda},\mu}$ of these differential operators by:
\begin{equation}\label{Lieder2}
X_h\cdot A:= L_{X_h}^{{\lambda},\mu}(A)=L_{X_h}^\mu\circ A-A\circ L_{X_h}^{{\lambda}}
\end{equation}
where $L_{X_h}^{{\lambda}}$ is the Lie derivative on $\mathcal{F}_{\lambda_1}\otimes\cdots\otimes\mathcal{F}_{\lambda_n}$  defined by the Leibniz rule.

Since the Lie algebra $\mathfrak{sl}(2)$ is a finite dimensional subalgebra of $\mathrm{Vect}(\mathbb{R})$, explicitely:
$$
\mathfrak{ sl}(2)=\mathrm{Span}(X_1,\,X_x,\,X_{x^2}),
$$
the spaces $\mathcal{F}_\mu$ and $\mathcal{D}_{{\lambda},\mu}$ are $\mathfrak{sl}(2)$-modules. 

\

Many authors were interested by the computation of the differential cohomology spaces $\mathrm{H}^j_\mathrm{diff}(\mathfrak{sl}(2), \mathcal{D}_{{\lambda},\mu})$ ($j=0,1,2$).
\begin{itemize}
\item[] If $n=1$, let us recall the computation by Feigin and Fuchs \cite{ff} of $\mathrm{H}^1_\mathrm{diff}(\mathrm{Vect}(\mathbb{R}), \Dlm)$, the computation by Gargoubi \cite{g} and Lecomte \cite{lec} of the spaces $\mathrm{H}^1_\mathrm{diff}(\mathfrak{sl}(2), \Dlm)$, while the $\mathfrak{sl}(2)$-relative cohomology spaces $\mathrm{H}^1_\mathrm{diff}(\mathrm{Vect}(\mathbb{R}),\mathfrak{sl}(2), \Dlm)$ were computed by Bouarroudj and Ovsienko \cite{bo}.\\

\item[] If $n=2$, the spaces $\mathrm{H}^1_\mathrm{diff}(\mathfrak{sl}(2),\Dlm)$ were computed by Bouarroudj \cite{b}.\\

\item[] If $n=3$, the spaces $\mathrm{H}^1_\mathrm{diff}(\mathfrak{sl}(2),\Dlm)$ were studied by O. Basdouri and E. Nasri \cite{be}.\\

\item[] For any $n$, the spaces $\mathrm{H}^1_\mathrm{diff}(\mathfrak{sl}(2), \Dlm)$ were studied by M. Ben Ammar and R. Sidaoui \cite{bs}.\\
\end{itemize}

Observe now that the space $M$ of $n$-differential operators with polynomial coefficients is a sub $\mathfrak{sl}(2)$-module of $\Dlm$ (a Harish-Chandra module). In the first section, we shall prove that the cohomologies $\mathrm{H}^j_{\mathrm{diff}}(\mathfrak{sl}(2),\Dlm)$ and $\mathrm{H}^j(\mathfrak{sl}(2),M)$ coincide. It is well known that the cohomology of this last Harish-Chandra module is localized in a sum $M^0+M^{-1}$ of weight spaces. More precisely, since the module $M$ is a `good' $\mathfrak{sl}(2)$-module in the sense of \cite{DMB}, the cohomology are characterized by an unique linear map $F:V\rightarrow V'$, where $V\subset M^0$ and $V'\subset M^{-1}$ are two explicit finite dimensional spaces, and:
$$\aligned
\HM0&\simeq\ker F,\\
\HM1&\simeq\ker F\oplus \coker F,\\
\HM2&\simeq\coker F.
\endaligned
$$

For completness, we recall the proof in Section \ref{map F}. Therefore, the computation of the dimension of one of these spaces gives immediately the dimension of the two others. 
Here we shall consider the dimension of $\HM2$. In a first step, we prove that these cohomology spaces are generically trivial. More precisely, we prove that the spaces $\mathrm{H}^0_{\mathrm{diff}}(\mathfrak{sl}(2),\Dlm)$ and $\mathrm{H}^1_{\mathrm{diff}}(\mathfrak{sl}(2),\Dlm)$ can be nontrivial only if $\mu-|\l|$ is in $\N$ while the space $\mathrm{H}^2_{\mathrm{diff}}(\mathfrak{sl}(2),\Dlm)$ can be nontrivial only if:
\begin{equation}\label{singular}
k=\mu-|{\lambda}|\in\mathbb{N}, \,\,-2{\lambda}\in\{0,\,\dots,\,k-1\}^n \,\,\text{and}\,\, -2|\lambda|\geq k-1.
\end{equation}
In Section \ref{mapf}, putting $-2\l=t=(t_1,\ldots,t_n)$, we restrict the map $F$ to a map $f$ from a vector subspace $W=W_{nkt}$ of $V$ to the space $W'=W_{n(k-1)t}\subset V'$, for which $\HM2=\coker f$. We then give an expression for the dimension $d_{nkt}$ of $W_{nkt}$ and prove, in Section \ref{maximalrank} that the map $f$ has always maximal rank. =Thus, for any $n$, $k$ and $t$,
$$
\dim \HM2=\sup(d_{n(k-1)t}-d_{nkt},0).
$$
We present a few examples.

In Section 6, we study the small $n$ cases, showing how computing this dimension for $n\leq 4$. Finally, in the last section, restricting ourselves to the so-called homogeneous case, we prove that the map $k\mapsto\dim\HM2$ is piecewise polynomial, with degree $n-2$.\\

%-----------------------------------------------------------------------------------------------------------------------

%-----------------------------------------------------------------------------------------------------------------------

\section{Various forms of cohomology }\label{various}

%-----------------------------------------------------------------------------------------------------------------------

Let $\mathfrak g$ be a Lie algebra and $V$ a $\mathfrak{g}$-module. For each $j=0,1,2,\dots$, we consider the space of $j$-cochains on $\mathfrak{g}$ with values in $V$:
\begin{equation*}
C^j(\mathfrak{g}, V ) := \mathrm{Hom}(\wedge^j(\mathfrak{g}),V)
\end{equation*}
and the coboundary operator $\partial=\partial_j: C^j(\mathfrak{g}, V)\rightarrow C^{j+1}(\mathfrak{g}, V )$ defined by 
$$\aligned
&(\partial c)(u_0,\dots,u_j)=\\
&=\sum_{i=0}^j (-1)^iu_i c(u_0\dots\hat{\imath}\dots u_j)+%\\
%&\hskip 2cm+ 
\hskip-0.2cm\sum_{0\leq i<i'\leq j}(-1)^{i+i'}c([u_i,u_{i'}],u_0\dots\hat{\imath}\dots\hat{\imath'}\dots u_j).
\endaligned
$$
Then $\partial^2=0$, and if $Z^j(\mathfrak{g},V)$ is the kernel of $\partial_j$ (the space of $j$ cocycles), it contains the image $B^j(\mathfrak{g},V)$ of $\partial_{j-1}$ (the space of $j$ coboundaries). As usual, the $j^{\mathrm{th}}$ space of cohomology of $\mathfrak g$ in $V$ is:
\begin{equation*}
\mathrm{H}^j (\mathfrak{g},V)=Z^j(\mathfrak{g},V)/B^j(\mathfrak{g},V).
\end{equation*}

In this paper, we study the cohomology of the space $\mathcal{D}_{{\lambda},\mu}$ of multi-differential operators (see the introduction) spanned, as a $C^\infty(\mathbb{R})$-module,  by the operators $\Omega^\alpha$, $\alpha=(\alpha_1,\dots,\alpha_n)\in\mathbb{N}^n$, defined by
$$
\Omega^\alpha(f_1dx^{\lambda_1}\otimes\cdots\otimes f_ndx^{\lambda_n})={1\over\alpha!}f_1^{(\alpha_1)}\dots f_n^{(\alpha_n)}dx^{\mu}.
$$
Recall that the Lie algebra $\mathfrak{sl}(2)$ is a subalgebra of $\mathrm{Vect}(\R)$, it is thus acting on this space by \eqref{Lieder2}
$$%\begin{equation}%\label{Lieder2}
X_h\cdot A:= L_{X_h}^{{\lambda},\mu}(A)=L_{X_h}^\mu\circ A-A\circ L_{X_h}^{{\lambda}}.
$$%\end{equation} 
Or, using the notation $t_i=-2\l_i$, $k=\mu-|\l|=\mu-\sum_i\l_i$,
\begin{equation}\label{action}
X_h\cdot A_\alpha\Omega^\alpha 
=  A'_\alpha h \Omega^{\alpha}+ (k-|\alpha|)A_\alpha h' \Omega^{\alpha}
-{1\over2}\displaystyle\sum_{i, \a_i>0}(\alpha_i-1-t_i) A_{\alpha}h''\Omega^{\alpha-\varepsilon_i},
 \end{equation}
where $(\varepsilon_1,\,\dots,\,\varepsilon_n)$ is the canonical basis of $\mathbb{R}^n$. 
We shall study the cohomology of this $\mathfrak{sl}(2)$-module. 

Before to begin, let us recall that some more restricted notion of cohomology were considered in the litterature. First, if we restrict ourselves to $j$-cochains (cocycles and coboundaries) which are the restrictions of $j$-differential operators on $\wedge^jC^\infty(\R)$ to  $\wedge^j\mathfrak{sl}(2)$, we define similarly the differential cohomology of $\mathcal{D}_{{\lambda},\mu}$, collection of the spaces denoted $\mathrm{H}^j_{\mathrm{diff}} (\mathfrak{sl}(2),\mathcal{D}_{{\lambda},\mu})$. This notion was used for instance in \cite{be,b}. 

It is also possible to restrict the $\mathfrak{sl}(2)$-action to the submodule $M$ of $n$-differential operators with polynomial coefficients:
$$
M=\{A=\sum_\alpha P_\alpha(x)\Omega^\alpha:~P_\alpha\text{ a polynomial function}\}.
$$
By construction, $M$ is a submodule of $\Dlm$. We can restrict ourselves to cochains with value in $M$ as in \cite{DMB}, introducing the corresponding cohomolgy. 

Moreover, we say that such a $j$-cochain $c:\wedge^j\mathfrak{sl}(2)\rightarrow M$ is reduced if, for any $u_2,\dots,u_j$ in $\mathfrak{sl}(2)$, $c(X_1,u_2,\dots,u_j)=0$. Let us denote by $\Zred{j}$ (resp. $\Bred{j}$) the space of reduced $j$-cocycles ($j$-coboundaries).

In the present section, we prove that these cohomologies are identical.\\

\begin{lem}\label{cocycles} 
For each $j$, we have :
$$\aligned
\ZD{j}&=\Zred{j}+\BD{j}\\
\ZDif{j}&=\Zred{j}+\BDif{j}\\
\ZM{j}&=\Zred{j}+\BM{j}.
\endaligned
$$
\end{lem}

\begin{proof}
%First, observe that the action is described by:
%\begin{equation}\label{action}X_h\cdot A_\alpha\Omega^\alpha 
%=  A'_\alpha h \Omega^{\alpha}+ (k-|\alpha|)A_\alpha h' \Omega^{\alpha}
%-{1\over2}\displaystyle\sum_{i=1}^n(\alpha_i+2\lambda_i-1) A_{\alpha}h''\Omega^{\alpha-\varepsilon_i}, \end{equation}
%where $(\varepsilon_1,\,\dots,\,\varepsilon_n)$ is the canonical basis of $\mathbb{R}^n$. 

Observe that, by \eqref{action}, the operator $X_1\cdot$ is simply the derivation operator:
$$
X_1\cdot A_\alpha(x)\Omega^\alpha=A'_\alpha(x)\Omega^\alpha.
$$
Therefore this operator is a surjective map and the space $M$ is:
$$
M=\{A\in\mathcal D_{{\lambda},\mu}:~(X_1)^q\cdot A=0\text{ for some natural number $q$}\}.
$$ 
Let us put $K=\ker(X_1\cdot)=\mathrm{Span}_\R\{\Omega^\alpha\}$, and observe that $X_x\cdot$ is diagolanized on $K$, since:
$$
X_x\cdot\Omega^\alpha=(k-|\alpha|)\Omega^\alpha.
$$
Denote by $K^m=\mathrm{Span}_\R\{\Omega^\alpha:~|\alpha|=m\}$ the eigenspace of $X_x\cdot$ in $K$.\\% with eigenvalue $k-m$.\\

\noindent
{$\mathbf{Case~} j=0$} 

There is nothing to prove if $j=0$. Indeed there is no coboundaries, and any $0$-cocycle $c\in\Dlm$ is differential and satisfies $X_1\cdot c=0$, $c$ is in $M$ and:
$$
\ZD{0}=\ZDif{0}=\ZM{0}=\Zred{0}.
$$

\noindent
{$\mathbf{Case~} j=1$} 

If $u\mapsto c(u)$ is 1-cochain, there is $b\in \Dlm$ such that $X_1\cdot b=c(X_1)$, if $c$ is a cocycle, the cocycle $c'=c-\partial b$ is thus reduced% (vanishing on $X_1$)
, the cocycle relation on $X_1$ and $X_x$ reduces to:
$$
0=\partial c'(X_1,X_x)=X_1\cdot c'(X_x)
$$
$c'(X_x)$ belongs to $K$, $c'(X_x)=\sum c'_\alpha\Omega^\alpha$. Similarly, on $X_1$ and $X_{x^2}$, we get: 
$$
0=\partial c'(X_1,X_{x^2})=X_1\cdot c'(X_{x^2})-c'([X_1,X_{x^2}])=X_1\cdot c'(X_{x^2})-2c'(X_x)
$$
this implies $X_1^2\cdot c'(X_{x^2})=0$, $c'(X_{x^2})$ is in $M$. Observe that  $b$ is differential. Moreover, if $c$ has values in $M$, we can choose $b\in M$, and we get:
$$\aligned
\ZD{1}&=\Zred{1}+\BD{1},\\
\ZDif{1}&=\Zred{1}+\BDif{1},\\
\ZM{1}&=\Zred{1}+\BM{1}.
\endaligned
$$

\noindent
{$\mathbf{Case~} j=2$} 

Suppose now $c$ is a $2$-cocycle, $c\in\ZD{2}$. We modify $c$ by a coboundary $\partial b$ where the cochain $b$ is reduced ($b(X_1)=0$) and satisfies $X_1\cdot b(X_x)=c(X_1,X_x)$ and $X_1\cdot b(X_{x^2})=c(X_1,X_{x^2})+2b(X_x)$. Then the $2$-cocycle $c'=c-\partial b$ satisfies:
$$\aligned
c'(X_1,X_x)&=c(X_1,X_x)-X_1\cdot b(X_x)=0,\\
c'(X_1,X_{x^2})&=c(X_1,X_{x^2})-X_1\cdot b(X_{x^2})-b([X_1,X_{x^2}])=0.
\endaligned
$$
Moreover the relation $\partial c'(X_1,X_x,X_{x^2})=0$ reduces to: $X_1\cdot c'(X_x,X_{x^2})=0$, this proves that $c'$ is reduced. 

It is also differential since, if $c'(X_x,X_{x^2})=\sum c'_\alpha D^\alpha$ (with $c'_\alpha$ constant), and if $h_i$ is a polynomial function with degree at most 2,
$$
c'(X_{h_1},X_{h_2})=\frac{1}{2}\left|\begin{matrix}h'_1&h'_2\\ h''_1&h''_2\end{matrix}\right|%(h'_1h''_2-h'_2h''_1)
\sum c'_\alpha D^\alpha.
$$
Observe that if $c$ had values in $M$, $b$ has also values in $M$, and if $c$ was differential, {\sl i. e.} if
$$\aligned
c(X_{h_1},X_{h_2})&=\left|\begin{matrix}h_1&h_2\\ h'_1&h'_2\end{matrix}\right|%(h_1h'_2-h_2h'_1)
A+\left|\begin{matrix}h_1&h_2\\ h''_1&h''_2\end{matrix}\right|%(h_1h''_2-h_2h''_1)
B+\left|\begin{matrix}h'_1&h'_2\\ h''_1&h''_2\end{matrix}\right|%(h'_1h''_2-h'_2h''_1)
C\\
&=\sum_\alpha\Big[\left|\begin{matrix}h_1&h_2\\ h'_1&h'_2\end{matrix}\right|%(h_1h'_2-h_2h'_1)
 A_\alpha+\left|\begin{matrix}h_1&h_2\\ h''_1&h''_2\end{matrix}\right|%(h_1h''_2-h_2h''_1) 
B_\alpha+\left|\begin{matrix}h'_1&h'_2\\ h''_1&h''_2\end{matrix}\right|%(h'_1h''_2-h'_2h''_1) 
C_{\alpha}\Big] \Omega^\alpha,
\endaligned
$$
then we choose $b$ also differential, and vanishing on $X_1$:
$$
b(X_h)=h'E+h''F=\sum_\alpha(h'E_\alpha+h''F_\alpha)\Omega^\alpha.
$$
A direct computation gives:  %{\it Sommes nous sûr de ce calcul ? les signes m'interrogent}
$$\aligned
\partial b (X_{h_1},X_{h_2})&=X_{h_1}\cdot b(X_{h_2})-X_{h_2}\cdot b(X_{h_1})-b([X_{h_1},X_{h_2}])\\
&=\sum_\alpha\Big[E'_\alpha \left|\begin{matrix}h_1&h_2\\ h'_1&h'_2\end{matrix}\right|%(h_1 h'_2-h_2 h'_1)
+F'_\alpha \left|\begin{matrix}h_1&h_2\\ h''_1&h''_2\end{matrix}\right|%(h_1h''_2-h_2h''_1)
+ \\
&\hskip-2cm+\big(\frac{1}{2} \sum_i(\alpha_i-t_i)E_{\alpha+\varepsilon_i}+(k-|\alpha|+1)F_{\alpha}\big)\left|\begin{matrix}h'_1&h'_2\\ h''_1&h''_2\end{matrix}\right|%(h'_1h''_2-h''_1h'_2)
\Big] \Omega^{\alpha}.
 \endaligned
$$
Choosing $E'_\alpha=A_\alpha$ and $F'_\alpha=B_\alpha$, we get a reduced, differential 2-cocycle $c'=c-\partial b$. We proved:
$$\aligned
\ZD{2}&=\Zred{2}+\BD{2},\\
\ZDif{2}&=\Zred{2}+\BDif{2},\\
\ZM{2}&=\Zred{2}+\BM{2}.
\endaligned
$$

\noindent
{$\mathbf{Case~} j=3$} 

If $c$  is a 3-cocycle, and $b$ is a 2-cochain such that $b(X_1,X_x)=b(X_1,X_{x^2})=0$ and $X_1\cdot b(X_x,X_{x^2})=c(X_1,X_x,X_{x^2})$, then $c=\partial b$. If $c$ has value in $M$, $b$ has also value in $M$, if $c$ is differential:
$$
c(X_{h_1},X_{h_2},X_{h_3})=\left|\begin{matrix}h_1&h_2&h_3\\ h'_1&h'_2&h'_3\\ h''_1&h''_2&h''_3\end{matrix}\right|\sum_\alpha E_\alpha\Omega^\alpha,
$$
then $b$ is also differential since 
$$
b(X_{h_1},X_{h_2})=\left|\begin{matrix}h'_1&h'_2\\ h''_1&h''_2\end{matrix}\right|\sum_\alpha F_\alpha\Omega^\alpha,
$$
with $F'_\alpha=E_\alpha$. Therefore,
$$\aligned
\ZD{3}&=\BD{3},\\
\ZDif{3}&=\BDif{3},\\
\ZM{3}&=\BM{3}.
\endaligned
$$
 
\noindent
{$\mathbf{Case~} j>3$} 

There is no $j$-cochain, and nothing to prove.\\
 \end{proof}

On the other hand,
\begin{lem}\label{coboundaries} 
For each $j$, we have :
$$\aligned
\BM{j}&=\ZM{j}\cap \BD{j},\\
%\BM{j}
&=\ZM{j}\cap \BDif{j}.\\
%\Bred{p}&=\Zred{p}\cap\BM{p}.
\endaligned
$$
\end{lem}

\begin{proof}
By definition, $\BM{j}\subset \ZM{j}\cap \BD{j}$ and $\BM{j}\subset \ZM{j}\cap \BDif{j}$. Let us prove the inverse inclusion. This is holding if $j\geq 3$ thanks to the proof of Lemma \ref{cocycles}, this holds also if $j=0$ since the only coboundary is then 0.\\

\noindent
{$\mathbf{Case~} j=1$} 

Let $c$ be in $\ZM{1}$ and $b$ in $\Dlm$ such that $c=\partial b$ (then $\partial b$ is in  $\BDif{1}=\BD{1}$). Then $X_1\cdot b=c(X_1)\in M$, therefore $b$ itself belongs to $M$, and $c$ in $\BM{1}$. This proves the lemma if $j=1$.\\

\noindent
{$\mathbf{Case~} j=2$} 

Let $c$ be in $\ZM{2}$ and $b$ a 1-cochain such that $c=\partial b$. There is $a$ in $\Dlm$ such that $b(X_1)=X_1\cdot a$, we replace $b$ by the reduced 1-cochain $b'=b-\partial a$, which is differential if $b$ was, we  get:
$$\aligned
c(X_1,X_x)&=X_1\cdot b'(X_x),\\
c(X_1,X_{x^2})&=X_1\cdot b'(X_{x^2})-2b(X_x).
\endaligned
$$
This implies that $b(X_1)$ and $b(X_x)$ belong to $M$, $b\in \BM{2}$, and proves the lemma if $j=2$.\\

\end{proof}

%Since $(Z+B)/B=Z/(Z\cap B)$, t
These two Lemmas give the identities between cohomologies:\\

\begin{cor}\label{cohomology}
For any $j$, we have:
$$
\HD{j}=\HDif{j}=\HM{j}.
$$
\end{cor}

\

%\begin{notation}
%From now on, we shall denote by $h^j_{nkt}$ the dimension of this space:
%$$
%h^j_{nkt}=\dim H^j(\mathfrak{sl}(2),M).
%$$
%\end{notation}

%-----------------------------------------------------------------------------------------------------------------------

%-----------------------------------------------------------------------------------------------------------------------

\section{The %cohomology and the 
map $F$}\label{map F}

%-----------------------------------------------------------------------------------------------------------------------

\ 

Recall that $K=\ker(X_1\cdot)=\Span\{\Omega^\alpha\}$ is the direct sum of the eigenspaces $K^m=\Span\{\Omega^\alpha:~|\alpha|=m\}$, for the eigenvalue $k-|\alpha|$ of $X_x\cdot$:
$$
X_x\cdot \Omega^\alpha=(k-|\alpha|)\Omega^\alpha.
$$ 

Then $u$ belongs to $K^k$ if and only if $X_1\cdot u=X_x\cdot u=0$. For such an $u$, consider $X_{x^2}\cdot u$, since:
$$
X_1\cdot(X_{x^2}\cdot u)= [X_1,X_{x^2}]\cdot u=2X_x\cdot u=0,
$$
$X_{x^2}\cdot u$ belongs to $K$. Similarly:
$$
X_x\cdot(X_{x^2}\cdot u)=[X_x,X_{x^2}]\cdot u=X_{x^2}\cdot u,
$$
then $X_{x^2}\cdot u$ belongs to $K^{k-1}$. We thus define the linear map $F:K^k\rightarrow K^{k-1}$ by $F(u)=X_{x^2}\cdot u$. The cohomology $\HD{*}$ is characterized by this map $F$. More precisely if $\coker(F)$ is the quotient $K^{k-1}/F(K^k)$,\\

\begin{lem}\label{fcohomology}

The spaces $\HD{j}$ are vanishing, except perhaps for:
$$\aligned
\HD{0}&=\ker(F),\\
\HD{1}&\simeq\ker(F)\oplus\coker(F),\\
\HD{2}&\simeq\coker(F).
\endaligned
$$
\end{lem}

\begin{proof}
We saw that $\HD{j}=0$ except perhaps for $j=0,1,2$.

\noindent
{$\mathbf{Case~} j=0$}

This is immediate, $c$ is in $\HD{0}$ if and only if $X_1\cdot c=X_x\cdot c=X_{x^2}\cdot c=0$, that is if and only if $c$ is in $\ker(F)$.

\noindent
{$\mathbf{Case~} j=1$}

Let $c:\mathfrak{sl}(2)\rightarrow \Dlm$ be a 1-cocycle. We saw that, up to a coboundary, we can choose $c$ reduced, {\sl i.e.} $c(X_h)\in M$ and $c(X_1)=0$, then the cocycle relation, $(\partial c)(X_1,X_x)=X_1\cdot c(X_x)$ gives $c(X_x)\in K$.
%$$
%(\partial c)(X_1,X_x)=X_1\cdot c(X_x)=0.
%$$
Thus $c(X_x)=\sum A_\alpha \Omega^\alpha$ ($A_\alpha$ constant), we can still add a coboundary $\partial b$ to $c$, such that $\partial b(X_1)=0$, by putting
$$
b=\sum_{|\alpha|\neq k}\frac{A_\alpha}{k-|\alpha|}\Omega^\alpha,
$$
then $c'=c-\partial b$ satisfies $c'(X_1)=0$, and $c'(X_x)\in K^k$. The cocycle relations are now:
$$\aligned
0=\partial c'(X_1,X_{x^2})&=X_1\cdot c'(X_{x^2})-c'([X_1,X_{x^2}])=X_1\cdot c'(X_{x^2})-2c'(X_x),\\
0=\partial c'(X_x,X_{x^2})&=X_x\cdot c'(X_{x^2})-X_{x^2}\cdot c'(X_x)-c'([X_x,X_{x^2}])\\
&=(X_x-1)\cdot c'(X_{x^2})-F(c'(X_x)).
\endaligned
$$
Since $F(c'(X_x))$ is an eigenvector for $X_x\cdot$, with eigenvalue 1, the cocycle relations are equivalent to 
$$\aligned
X_1\cdot c'(X_{x^2})&=2c'(X_x),\\
(X_x-1)\cdot c'(X_{x^2})&=0,\\
F(c'(X_x))&=0.\\
\endaligned
$$
The two first relations give $[c'(X_{x^2})-2xc'(X_x)]$ is in $K^{k-1}$. We can still add to $c'$ a coboundary $\partial b'$ such that $X_1\cdot b'=X_x\cdot b'=0$, ({\sl i.e.} $b'\in K^k$), then $c''=c'-\partial b'$ satisfies $c''(X_1)=c'(X_1)=0$, $c''(X_x)=c'(X_x)\in\ker(F)$ and $[c''(X_{x^2})-2xc''(X_x)]=[c'(X_{x^2}))-2xc'(X_x)]-F(b')$. Identifying $\coker(F)$ to a fixed subspace in $K^{k-1}$, supplementary to $F(K^k)$, we proved that, up to a coboundary, we can choose the cocycle $c$ such that
$$
c(X_1)=0,~~c(X_x)\in\ker(F),~~c(X_{x^2})-2xc(X_x)\in\coker(F).
$$
Now, a direct computation proves that any map $c:\mathfrak{sl}(2)\rightarrow \Dlm$ satisfying these relations defines a 1-cocycle, and $c$ is a coboundary if and only if $c(X_x)=c(X_{x^2})=0$. This proves
$$
\HD{1}=\ker(F)\oplus\coker(F).
$$

\noindent
{$\mathbf{Case~} j=2$}

If $c$ is a 2-cocycle, we can suppose that $c$ is reduced ($c(X_1,X_x)=c(X_1,X_{x^2})=0$), then $c(X_x,X_{x^2})=\sum_\alpha A_\alpha \Omega^\alpha$ is in $K$;
%$$
%c(X_x,X_{x^2})=\sum_\alpha A_\alpha \Omega^\alpha
%$$
we can add to $c$ a coboundary $\partial b$ with $b(X_1)=b(X_x)=0$, and
$$
b(X_{x^2})=\sum_{|\alpha|\neq k-1}\frac{A_\alpha}{|\alpha|-k+1}\Omega^\alpha.
$$
Then $c'=c-\partial b$ satisfies now:
$$\aligned
c'(X_1,X_x)&=0,\\
c'(X_1,X_{x^2})&=0,\\
c'(X_x,X_{x^2})&\in K^{k-1}.
\endaligned
$$
Writing as above $K^{k-1}=F(K^k)\oplus \coker(F)$, we can finally add to $c'$ a coboundary $b'$ such that $b'(X_1)=b'(X_x)=0$, $b'(X_{x^2})\in K^k$, such that $c''=c'-\partial b'$ satisfies $c''(X_1,X_x)=c''(X_1,X_{x^2})=0$, and $c''(X_x,X_{x^2})\in \coker(F)$.

Now a direct computation shows that any bilinear skew-symmetric map $c:\mathfrak{sl}(2)\wedge\mathfrak{sl}(2)\rightarrow \Dlm$ such that $c(X_1,X_x)=c(X_1,X_{x^2})=0$, and $c(X_x,X_{x^2})\in \coker(F)$, is a 2-cocycle, moreover if $c$ is a coboundary $\partial b$, it is a coboundary $\partial b'$ with $b'(X_1)=b'(X_x)=0$ and $b'(X_{x^2})\in K^k$, then $c=0$. This achieves the proof of the lemma.\\
\end{proof}

%-----------------------------------------------------------------------------------------------------------------------

%-----------------------------------------------------------------------------------------------------------------------
\section{General case}

\
%-----------------------------------------------------------------------------------------------------------------------

\begin{notation}
From now on, we change our notation by putting:
$$
k=\mu-|\l|,\hskip 0.7cm t_i=-2\l_i,~~t=(t_1,\ldots,t_n).
$$

%since $K^k\neq 0$ only if $k$ is a non negative integral number, we suppose $k$ is in $\N$.\\

We now use these parameters everywhere, for instance we put:
$$
h^j=h^j_{nkt}=\dim\HD{j}.% H^j(\mathfrak{sl}(2),\Dlm)~~(j=0,1,2).
$$
\end{notation}

Using the relation \eqref{action}, we see that, if $k$ is a positive integer, a basis of $K^k$ is $\{\Omega^\alpha:~|\alpha|=k\}$, and $F$ is given by:
\begin{equation}\label{F}
F(\Omega^\alpha)=\sum_{i:~\a_i\neq0} (t_i+1-\a_i)\Omega^{\alpha-\varepsilon_i}.
\end{equation}
On the other hand, if $k$ is not in $\N$, $K^k=0$, and $F$ is the zero map.\\

Suppose that $k\in\N$, then $\dim K^k=\#\{\alpha:~|\alpha|=k\}=\left(\begin{smallmatrix}n+k-1\\ k\end{smallmatrix}\right)$, therefore 
$$\aligned
h^2&=\dim\coker(F)\\
&=\left(\begin{matrix}n+k-2\\ k-1\end{matrix}\right)-\dim F(K^k)\\
&=\left(\begin{matrix}n+k-2\\ k-1\end{matrix}\right)-\left(\begin{matrix}n+k-1\\ k\end{matrix}\right)+\dim\ker(F)\\
&=h^0-\left(\begin{matrix}n+k-2\\ k\end{matrix}\right).
\endaligned
$$
Recall that $h^1=h^0+h^2$. Hence, the computation of $h^2=h^2_{nkt}$ gives a complete description of the cohomology $\HD{*}$.\\

%Let us now prove that we can restrict ourselves to the case $k\in\N$ and $t_i\in\{0,\ldots,k\}$.\\

\begin{thm}\label{third} ~~

\begin{itemize}
\item [1.] If $%\mu-|\underline{\lambda|}
k\notin\N$ then $\HD{j}=0%\mathrm{H}^j_\mathrm{diff}(\mathfrak{sl}(2),\mathcal{D}_{\underline{\lambda},\mu})= 0
$ for any $ j$.

\item [2.] If $k\in\N$ and there is $i$ such that $t_i\notin\{0,\ldots,k-1\}$, then $F$ is surjective, $h^2_{nkt}=0$.

\item[3.] Suppose $k\in\N$ and $t_1\notin \{0,\ldots,k-1\}$, then the spaces $\HD{0}$ and $\HD{1}$ are both isomorphic to the subspace of $K^k$ spanned by the $\O^\a$ such that $\a_1=0$ and $|\a|=k$, their common dimension is $\left(\begin{smallmatrix}n+k-2\\ k\end{smallmatrix}\right)$.\\
\end{itemize}
\end{thm}

\begin{proof} Recall that 
$$
\HD{2}\simeq\coker(F)
$$
where $F$ is the linear map $F:K^k\rightarrow K^{k-1}$ defined in \eqref{F}.
\begin{itemize}
\item[1.] Therefore, the space $\HDif{2}$ is trivial if and only if the linear map $F$ is surjective. Of course, if $k\notin\N$ then $K^{k-1}=0$ and $F$ is onto. 

\item[2.] Suppose now $k\in\N$  and $t_i\notin\{0,\ldots,k-1\}$ for some $i$. By modifying the ordering of the indexes $i$, we can suppose $i=1$. We put thus the lexicographic ordering on the multi-indexes $\a$ such that $|\a|=k-1$. Then the largest such index is $(k-1)\varepsilon_1$, and the corresponding vector $\O^{(k-1)\varepsilon_1}$ is in $F(K^k)$ since :
$$
F(\O^{k\varepsilon_1})=(t_1+1-k)\O^{(k-1)\varepsilon_1}\qquad\text{ and }\quad t_1+k-1\neq0.
$$
Now suppose there is $\beta$ such that $|\beta|=k-1$ and $\O^\a$ is in $F(K^k)$ for any $\a>\beta$. Then we have:
$$
F(\O^{\beta+\varepsilon_1})=(t_1-\beta_1)\O^\beta+\sum_{i>1,~\beta_i\neq0}(\beta_i+1-t_i)\O^{\beta+\varepsilon_1-\varepsilon_i}
$$
Since $\beta_1\in\{0,\ldots,k-1\}$ and $\beta+\varepsilon_1-\varepsilon_i>\beta$, this proves that $\O^\beta$ is in $F(K^k)$. Therefore any $\O^\a$ is in $F(K^k)$, the map $F$ is surjective.

\item[3.] Suppose now $k\in\N$ and $t_1\notin\{0,\ldots,k-1\}$, then $F$ is onto, thus $\HD{0}$ and $\HD{1}$ are both isomorphic to $\ker(F)$, that is to the space of solution of the linear system in the variables $A_\a$ ($|\a|=k$):
$$
0=F(\sum_{|\a|=k} A_\a\O^\a)=\sum_{|\beta|=k-1}\big(\sum_i(\beta_i-t_i)A_{\beta+\varepsilon_i}\big)\O^{\beta}.
$$
In fact, this system is triangular by block, it is the union of the sub-system defined for each $0\leq b\leq k-1$, by the equations for $\beta_1=b$, which are, if we write $\beta=(b,\beta')$ and similarly $\a=(a,\a')$:
$$
0=(b-t_1)A_{(b+1,\beta')}+\sum_{i>1}(\beta'_i-t_i)A_{(b,\beta'+\varepsilon_i)}\qquad (|\beta'|=k-1-b).
$$ 
The solutions of our system are therefore parametrized by the vectors $Z=\sum_{|\alpha'|=k}z_{(0,\alpha')}\O^{(0,\alpha')}$, and they are recursively defined for each $a$ by:
$$
\left\{\aligned A_{(0,\alpha')}&=z_{(0,\alpha')}&\hskip1cm &(|\alpha'|=k)\\
A_{(a+1,\alpha')}&=-\sum_{i>1}\frac{\alpha'_i-t_i}{a-t_1}A_{(a,\a'+\varepsilon_i)} &&(|\a'|=k-1-a).
\endaligned\right.
$$
Denote by $S(Z)$ this solution. The map $S$ is a linear isomorphism between the space $K^{0,k}$ with basis $(\O^{(0,\alpha')})_{|\alpha'|=k}$, whose dimension is $\left(\begin{smallmatrix}n+k-2\\ k\end{smallmatrix}\right)$ and $\ker(F)$. This ends the proof of the theorem.\\
\end{itemize}
\end{proof}

More precisely, using Lemma \ref{fcohomology}, we get:

\begin{cor}\label{generic cocycles}

Suppose $k\in\N$, and $t_1\notin\{0,\ldots,k-1\}$ then the space $\HD{0}$ is the space of operators $S(Z)$,
%\begin{equation*}\label{cocycth}
%A=\sum_{|\alpha|=k} S(Z)_\alpha  \Omega^{\alpha},
%\end{equation*}
where $Z=\sum_{|\alpha'|=k}z_{(0,\alpha')}\O^{(0,\alpha')}$ is any vector in $K^{0,k}$.

Similarily, the space $\HD{1}$ is the space of  the applications $c_Z$ defined by $
%\begin{equation*}\label{cocycth}
c_Z(X_h) =h'S(Z)$,
%\end{equation*}
where $Z$ is any vector in $K^{0,k}$.\\
\end{cor}

%\begin{proof} 
%The proof follows from Lemma \ref{fcohomology}. Note that in this case the map $F$ is surjective, therefore, we can construct the cocycles $c$ generating $\HD{1}$   by choosing $c(X_1)=0$, $c(X_x)\in \ker F$ and $c(X_{x^2})=2xc(X_x)$.\\
%\end{proof}

From now on, we suppose $k\in\N$, $t_i\in\N$ and $k>t_i$ for each $i$.\\

%-----------------------------------------------------------------------------------------------------------------------

%-----------------------------------------------------------------------------------------------------------------------

\section{A restriction $f$ of $F$}\label{mapf}

\

%-----------------------------------------------------------------------------------------------------------------------

Suppose as above $k\in\N$, $t_i\in\N$ and $k>t_i$ for each $i$.  Let us now restrict our study to a sub-basis of $\{\Omega^\alpha\}$ and the corresponding restriction $f$ of $F$.\\

\begin{lem}
Recall the map $F:K^k\rightarrow K^{k-1}$, describing the cohomology $\HD{*}$, and the basis $\mathcal B^m=\{\Omega^\alpha:~|\alpha|=m\}$ of $K^m$ ($m=k$ or $k-1$). Denote
$$
\mathcal B^m_t=\{\Omega^\alpha\in\mathcal B^m:~\alpha_i\leq t_i ~(i=1,\ldots,n)\},
$$
and let $K_t^m$ be the vector space generated by $\mathcal B^m_t$. Then:
$$
F(K_t^k)\subset K_t^{k-1}\quad\text{ and }\quad K^{k-1}=K_t^{k-1}+F(K^k).
$$
Denote by $f:K_t^k\rightarrow K_t^{k-1}$ the restriction of the map $-F$. Then:
$$
\coker(f)\simeq \coker(F).
$$
\end{lem}

\begin{proof} As in the proof of Proposition \ref{third}, we put on the multi-indexes $\alpha$ (and the corresponding vectors $\Omega^\alpha$) the inverse lexicographic ordering. Then, the smallest element in $\mathcal B^{k-1}$ is $\Omega^{(k-1)\varepsilon_1}$, if $k-1>t_1$, it is in $F(K^k)$. Suppose now there is $a>t_1$ such that each $\Omega^\alpha\in\mathcal B^{k-1}$, with $\alpha_1\geq a$ is in $F(K^k)$. Let us consider the smallest $\Omega^\alpha$, with $\alpha_1<a$ in $\mathcal B^k$. Then:
$$
F(\Omega^{\a+\varepsilon_1})=(t_1-\alpha_1)\Omega^\alpha+\sum_{i=2}^n (t_i+1-\beta_i)\Omega^{\alpha+\varepsilon_1-\varepsilon_i}.
$$
Our induction hypothesis implies that each $\Omega^{\alpha+\varepsilon_1-\varepsilon_i}$ is in $F(K^k)$, therefore $\Omega^\alpha$ itself is $F(K^k)$. This proves that each $\Omega^\alpha$ in $\mathcal B^{k-1}$ such that $\alpha_1>t_1$ is in $F(K^k)$. 

Denote by $\mathcal B^{k-1}_{>t}$ the set of the $\Omega^\alpha$ such that $\alpha_i>t_i$ for some $i$, using a permutation, the preceding proves that $\mathcal B^{k-1}_{>t}\subset F(K^k)$, or, if $K^{k-1}_{>t}$ is the vector space spanned by $\mathcal B^{k-1}_{>t}$,
$$
K^{k-1}=K_t^{k-1}\oplus K^{k-1}_{>t}=K_t^{k-1}+F(K^k).
$$
Let now $\Omega^\alpha$ be a vector in $\mathcal B^k$, since:
$$
F(\Omega^\alpha)=\sum_{i: ~\a_i>0} (t_i+1-\alpha_i)\Omega^{\alpha-\varepsilon_i},
$$
it is clear that $F(K^k_t)\subset K^{k-1}_t$, moreover, if there is $j$ such that $\alpha_j>t_j+1$, then $F(\Omega^\alpha)$ is in $K^{k-1}_{>t}$, similarly, if there is $j$ such that $\alpha_j=t_j+1$, then
$$
F(\Omega^\alpha)=\sum_{i\neq j,~\a_i>0} (t_i+1-\alpha_i)\Omega^{\alpha-\varepsilon_i}\in K^{k-1}_{>t}.
$$
In other words, $F(K^k_{>t})\subset K^{k-1}_{>t}$, thus $F(K^k_t)\cap F(K^k_{>t})=0$ and:
$$
F(K^k)\cap K^{k-1}_t=F(K^k_t),\quad\quad
%$$
%In other words:
%$$
\coker (F)%=K^{k-1}_t/(F(K^k)\cap K^{k-1}_t)
=\coker(f).
$$
\end{proof}

In the next section, we shall prove that the rank of $f$ is always maximal. Just now, let us present the matrix $A_{nkt}$ of $f$ in a simple form.\\

We first change the basis in $K^k_t$ and $K^{k-1}_t$, by putting:
\begin{equation}\label{Eformula}
E_\alpha=\prod_{i=1}^n\frac{1}{(\alpha_i-1-t_i)^{\alpha_i} } \Omega^\alpha
\end{equation}
Thus
$$\aligned
f(E_\alpha)%&=\sum_{j:~\alpha_j>0} \prod_{i=1}^n\frac{1}{(\alpha_i-1-t_i)^{\alpha_i}}(\alpha_j-1-t_j)\Omega^{\alpha-\varepsilon_j}\\
%&=\sum_{j:~\alpha_j>0}\prod_{i\neq j}\frac{1}{(\alpha_i-t_i-1)^{\alpha_i}}~\frac{1}{(\alpha_j-t_j-1)^{\alpha_j-1}}\Omega^{\alpha-\varepsilon_j}\\
&=\sum_{j:~\alpha_j>0} E_{\alpha-\varepsilon_j}.
\endaligned
$$

Let us denote by $A_{nkt}$ the matrix of $f$, on these basis. An entry $x_{\alpha',\alpha}$ in $A_{nkt}$ is indexed by $n$-indexes $\alpha$ (for the columns) and $\alpha'$ for the rows, with $0\leq\alpha_i,\alpha'_i\leq t_i$  for each $i$ and $|\alpha|=k$, $|\alpha'|=k-1$. Let $D_{nkt}$ (resp. $D_{n(k-1)t}$) these sets of multi indexes. Put:
$$
d_{nkt}=\#D_{nkt}=\#\{\alpha:~\alpha_i\leq t_i,~|\alpha|=k\}.
$$

The matrix $A_{nkt}$ is a $d_{n(k-1)t}\times d_{nkt}$-matrix, its entries $x_{\alpha',\alpha}$ are:
$$
x_{\alpha',\alpha}=\left\{ \aligned &1\text{ if there is }j\text{ such that }\alpha'=\alpha-\varepsilon_j\\
&0\text{ else.}\endaligned  \right.
$$

%Observe now that the map $\alpha\mapsto \alpha^*=\beta$ defined by $\beta_i=t-\alpha_i$ (thus $|\beta|=nt-|\alpha|$), with inverse map  $\alpha=\beta^*$ ($\alpha_i=t-\beta_i$for each $i$. Therefore $d_{nkt}=d(n,nt-k,t)$. By construction, the map $*$ is decreasing for the inverse lexicographic ordering: $\alpha< \beta$ is equivalent to $\beta^*<\alpha^*$.

For any $p$, let $J=J_p$ be the $p\times p$ matrix such that $J_{ij}=0$ except for the entries $J_{i(p+1-i)}$ ($1\leq i\leq p$), of course, $J^t=J$ and $J^2=I_p$ (the identity matrix). On the space $\R^p$, we define the symmetric bilinear form:
$$
\langle x,y\rangle=x^tJy.
$$
It is non degenerated. More precisely if $\Vert ~\Vert$ is the usual Euclidean norm,
$$
\langle Jx,x\rangle=x^tJ^2x=\Vert x\Vert^2.
$$
Now if $A$ is a $r\times s$ matrix, we define the $s\times r$ matrix $A^*$ by $A^*=J_sA^tJ_r$.

The usual relations hold, if $B$ is a $s\times t$ matrix, $(A^*)^*=A$, $(AB)^*=B^*A^*$.

Moreover, for any $x$ and $y$ in $\R^p$, viewed as $p\times 1$ matrices, $x^*y=\langle x,y\rangle$.
%$$
%\langle x,y\rangle=x^tJ_py=J_1x^tJ_py=x^*y.
%$$
Therefore, if $y\in\R^r$, and $x\in\R^s$,
$$
\langle x,Ay\rangle=x^*Ay=(A^*x)^*y=\langle A^*x,y\rangle.
$$
Let us say that a matrix $A$ is injective (resp. surjective) if the corresponding linear map is injective (resp. surjective) or if its columns  (resp. its rows) are independent. We shall use the following simple result:

\begin{lem}\label{injcritere}
Let $A$ be a $r\times s$ matrix, then $A$ is injective if and only if its adjoint $A^*$ is surjective. Moreover, if $A$ is injective, then the $r\times r$ matrix $A^*A$ is invertible.\\
\end{lem}

\begin{proof}
The first assertion is evident. Suppose now $A$ injective, observe that $A^*A$ is a square matrix, and let $x$ be in the kernel of $A^*A$, then:
$$
0=x^*A^*Ax=(Ax)^*(Ax)=\Vert Ax\Vert^2,
$$
therefore $Ax=0$, and $x=0$.\\
\end{proof}

Let us now compute the adjoint of the matrix $A_{nkt}$.%: it is a $d(n,k-1,t)\times d(n,k,t)$-matrix, its entries $x_{\alpha',\alpha}$ are:
%$$
%x_{\alpha',\alpha}=\left\{ \aligned &1\text{ if there is }j\text{ such that }\alpha'=\alpha-\varepsilon_j\\
%&0\text{ else.}\endaligned  \right.
%$$
Consider the map $\alpha\mapsto \alpha^*=\beta$ defined by $\beta_i=t_i-\alpha_i$ (thus $|\beta|=|t|-|\alpha|=|t|-k$), with inverse map  $\alpha=\beta^*$. %($\alpha_i=t_i-\beta_i$ for each $i$). 
Especially $d_{nkt}=d_{n(|t|-k)t}$. By construction, the map $*$ is decreasing for the inverse lexicographic orderings: $\alpha< \beta$ if and only if $\beta^*<\alpha^*$.

Fix $n$ and $t$, put $d_k=d_{nkt}$, we now define a $d_k\times d_{k-1}$-matrix $B=[y_{\beta',\beta}]$ by putting $y_{\beta',\beta}=x_{\beta^*,{\beta'}^*}$. Then 
$$
y_{\beta',\beta}=\left\{ \aligned &1\text{ if there is }j\text{ such that }\beta^*={\beta'}^*-\varepsilon_j\text{ or }\beta'=\beta-\varepsilon_j,\\
&0\text{ else.}\endaligned  \right.
$$
In other words, $B=A_{n(|t|+1-k)t}$, and since $y_{\beta',\beta}=x_{\beta^*,{\beta'}^*}$, and taking into account the inversion of ordering on the indexes, we also get $B^t=J_{d_{k-1}}A(n,k,t)J_{d_k}$, this proves:

\begin{lem}\label{adjointofA}
The matrix, adjoint of $A_{nkt}$ is:
$$
A_{nkt}^*=A_{n(|t|+1-k)t}.
$$
\end{lem}

Observe that the sets $D_{nkt}$ and $D_{n(k-1)t}$ are both empty if $|t|+1<k$. So we suppose now $1\leq k\leq |t|+1$ (the case $k=0$ being trivial). Using now the lemma \ref{adjointofA}, it is enough to study the rank of $A_{nkt}$ for $\frac{|t|+1}{2}\leq k\leq|t|+1$. The next section will be devoted to this point.\\

%-----------------------------------------------------------------------------------------------------------------------

\section{The rank of $f$ is maximal}\label{maximalrank}

\
%-----------------------------------------------------------------------------------------------------------------------

Thank to Lemma \ref{adjointofA}, to prove that the rank of $f$ is maximal, it is enough to prove that, for any $n\geq1$, any $t=(t_1,\ldots,t_n)$ and any $k$ such that $\frac{|t|+1}{2}\leq k\leq|t|+1$, the matrix $A_{nkt}$ is injective.

Then, computing the $d_{nkt}$, we get the dimension of the second cohomology space $\HD{2}$.\\

Recall that $d_{nkt}=\# D_{nkt}$ where
$$
D_{nkt}=\{\alpha=(\alpha_1,\ldots,\alpha_n)\in \N^n:~\alpha_i\leq t_i~(1\leq i\leq n),~|\alpha|=k\}.
$$
We shall use an induction on $n$. More precisely, we shall prove that $f_{nkt}$ has maximal rank, supposing this is holding for any $f_{n'k't'}$ with $n'=n-1$, $t'=(t_2,\dots,t_n)$, and $k'=k-1$, however we do not suppose $k>t_i$ here.\\

First, let us restrict the choice of the $t_i$, thanks to the following simple Lemma:

\begin{lem}\label{permutation}
For any permutation $\sigma$ of $\{1,\ldots,n\}$, we put $t^\sigma=(t_{\sigma(1)},\ldots,t_{\sigma(n)})$. Then there is an isomorphism $\varphi_{nkt}^\sigma:W_{nkt^\sigma}\rightarrow W_{nkt}$ such that (with $k'=k-1$)
$$
f_{nkt}=(\varphi^\sigma_{nk't^\sigma})^{-1}\circ f_{nkt^\sigma}\circ \varphi_{nkt}^\sigma.
$$

%On the other hand, if $t_1=0$ then $f_{nkt}\simeq f_{n'k't'}$.\\
\end{lem}

%\begin{proof}
%The map $\varphi^\sigma_{nkt}: E_\alpha\mapsto E_{\alpha^\sigma}$ defines a bijection from the basis $\{E_\alpha:~\alpha\in D_{nkt}\}$ of $W_{nkt}$ onto the basis $\{E_\beta:~\beta\in D_{nkt^\sigma}\}$ of $W_{nkt^\sigma}$. Therefore the spaces $W_{nkt}$ and $W_{nkt^\sigma}$ have the same dimension%: $d_{nkt}=d_{nkt^\sigma}$
%, and $\varphi^\sigma_{nkt}$ gives rise to an invertible linear map from $W_{nkt}$ onto $W_{nkt^\sigma}$. By construction $\varphi^{\sigma^{-1}}_{nkt^\sigma}=(\varphi^\sigma_{nkt})^{-1}$ and 
%$$
%f_{nkt}=(\varphi^\sigma_{nk't^\sigma})^{-1}\circ f_{nkt^\sigma}\circ \varphi_{nkt}^\sigma.
%$$
%Therefore $f_{nkt}$ has maximal rank if and only if $f_{nkt^\sigma}$ has maximal rank.

%Suppose now $t_1=0$. Then each multi-index $\alpha\in D_{nkt}$ has the form:
%$$
%\alpha=(0,\alpha_2,\ldots,\alpha_n)=(0,\beta),
%$$
%with $\beta\in D_{n'kt'}$. Thus, identifying $E_{(0,\beta)}\in W_{nkt}$ with $E_\beta\in W_{n'kt'}$, we get $d_{nkt}=d_{n'kt'}$. Observe now that, with this identification, $f_{nkt}(E_{(0,\beta)}=f_{n'kt'}(E_\beta)$. Therefore the rank of $f_{nkt}$ is maximal if and only if the rank of $f_{n'kt'}$ is maximal.\\ 

%\end{proof}

Using this lemma, in our induction proof, we can now suppose that
$$
t_1\leq t_2\leq\dots\leq t_n
$$ 
holds. Then we also have $n>1$ and $t_1\leq t_2+\dots+t_n=|t'|$.\\

Observe that $D_{nkt}$ is a disjoint union:
$$
D_{nkt}=\bigsqcup_{a=0}^{t_1}\{(a,\alpha_2,\ldots,\alpha_n):~\alpha_i\leq t_i,~\alpha_2+\dots+\alpha_n=k-a\}.
$$
(Some of these subsets can be empty). Therefore, with our notation:
$$
d_{nkt}=\sum_{a=0}^{t_1}d_{n'(k-a)t'}.
$$

Let us accordingly decompose the matrix $A_{nkt}$ into blocks  $X_{a'a}$ in which the indexes $\alpha$ (resp. $\alpha'$) satisfy $\alpha_1=a$ (resp. $\alpha'_1=a'$). Observe that, if $n'=n-1$ and $t'=(t_2,\ldots,t_n)$,
\begin{itemize}
\item[1.] If $a-a'\notin\{0,1\}$, then $X_{a'a}=0$. 

\item[2.] If $a'=a-1$, then the corresponding indexes $\alpha$ and $\alpha'$ can be written $\alpha=(a,\alpha_2,\ldots,\alpha_n)$ (resp. $\alpha'=(a-1,\alpha'_2,\ldots,\alpha'_n)$. Therefore $x_{\alpha',\alpha}$ does not vanish if and only if $(\alpha_2,\ldots,\alpha_n)=(\alpha'_2,\ldots,\alpha'_n)$. Since $\alpha_2+\dots+\alpha_n=k-a$, hence $X_{a-1,a}$ is the identity matrix with size $d_{n'(k-a)t'}$:
$$
X_{a-1,a}=I_{d_{n'(k-a)t'}}=:I_{(k-a)}.
$$

\item[3.] If $a=a'$, we write $\alpha=(a,\alpha_2,\ldots,\alpha_n)$ ($\alpha_2+\dots+\alpha_n=k-a$) and $\alpha'=(a,\alpha'_2,\ldots,\alpha'_n)$ ($\alpha'_2+\dots+\alpha'_n=k-a-1$). Looking only on the $n-1$ last indexes, we get:
$$
X_{a,a}=A_{n'(k-a)t'}=:A_{(k-a)}.
$$
\end{itemize}

\begin{thm}\label{rankf}
The rank of the linear map $f=f_{nkt}:W_{nkt}\rightarrow W_{n(k-1)t}$ is always maximal. More precisely, if $\frac{|t|+1}{2}\leq k\leq|t|+1$, then $f_{nkt}$ is injective, and if $1\leq k\leq \frac{|t|+1}{2}$, then $f_{nkt}$ is surjective.\\
\end{thm}

\begin{proof}
If $n=1$, and $t<k'=k-1$, then $f_{1kt}$ does not exist, if $t=k'$, $f_{1kk'}$ is the 0 map from $\{0\}$ into the 1-dimensional space $\C E_{\varepsilon_1}$, it is injective. If $t<k'$, $f_{1kt}=\mathrm{Id}$ is invertible. From now on, we suppose $n\geq 2$ and Theorem \ref{rankf} holds for $n-1$. We suppose also $\frac{|t|+1}{2}\leq k\leq|t|+1$.\\

Put $n'=n-1$ and $t'=(t_2,\ldots,t_n)$. Since:
$$
\frac{|t'|+t_1+1}{2}=\frac{|t|+1}{2}\leq k\leq |t|+1,
$$
then the inequality $k\leq t_1$ is not holding, since this would imply $k\leq \frac{|t'|+t_1}{2}<\frac{|t'|+t_1+1}{2}$. 

Looking to the relative positions of $k$, $t_1$ and $|t'|$, the only possibilities are 
%\begin{itemize}
%\item[(1)] $t_1\leq |t'|<k$,
%\item[(2)] $t_1<k\leq |t'|$, or $\frac{|t|+1}{2}\leq k\leq |t'|$.
%\end{itemize}

%*************
%Moreover, in the last case, we shall consider the two subcases :
%\begin{itemize} 
%\item[(2.1)] $\frac{|t'|+1}{2}\leq k-t_1$,
%\item[(2.2)] $k-t_1<\frac{|t'|+1}{2}$.
%\end{itemize}
%*************

\begin{itemize}
\item[\bf Case (1)] $t_1\leq|t'|<k$.

In this case, the matrix $A_{nkt}$ has the following form:
$$
A_{nkt}={\tiny\begin{tabular}{c|cccccc}
\diagbox{$a'$}{$a$}&$t_1$&$t_1-1$&$\cdots$&$k-|t'|+1$&$k-|t'|$\\
\hline
$t_1$&$A_{(k-t_1)}$&0&$\cdots$&$0$&$0$\\
$t_1-1$&$I_{(k-t_1)}$&$A_{(k-t_1+1)}$&&&\\
$\vdots$&&$\ddots$&$\ddots$&&\\
%$\vdots$&&&$\ddots$&$\ddots$&&&\\
%$a$&&&&$B_{d_{k-a-1},d_{k-a}}$&$0$&&&\\
%$a-1$&&&&$I_{k-a}$&$B_{k-a+1}$&&&\\
$\vdots$&&&$\ddots$&$\ddots$&\\
%$\vdots$&&&&&&$\ddots$&$\ddots$&\\
$k-|t'|$&$0$&$0$&$\cdots$&$I_{(|t'|-1)}$&$A_{(|t'|)}$\\
$k-|t'|-1$&$0$&$0$&$\cdots$&$0$&$I_{(|t'|)}$\\
\end{tabular}}.
$$
The matrix is in echelon form, its columns are independent, $A_{nkt}$ is injective. Moreover the dimension of its coker is the number of rows of $A_{(k-t_1)}$, we get ($k'=k-1$)
$$
h^2_{nkt}=\dim(\coker(f_{nkt}))=d_{nk't}-d_{nkt}=d_{n'(k'-t_1)t'}.
$$
But since here $D_{n'kt'}=\emptyset$, we also can write:
$$
h^2_{nkt}=d_{n'(k'-t_1)t'}-d_{n'kt'}.
$$

\

\item[\bf Case (2)] $t_1<k\leq |t'|$.

Now $t_1\leq k-1$ and the matrix $A_{nkt}$ is:
$$
A_{nkt}={\tiny\begin{tabular}{c|cccccc}
\diagbox{$a'$}{$a$}&$t_1$&$t_1-1$&$\cdots$&$1$&$0$\\
\hline
$t_1$&$A_{(k-t_1)}$&0&$\cdots$&$0$&$0$\\
$t_1-1$&$I_{(k-t_1)}$&$A_{(k-t_1+1)}$&&&\\
$\vdots$&&$\ddots$&$\ddots$&&\\
$\vdots$&&&$\ddots$&$\ddots$&\\
$0$&$0$&$0$&$\cdots$&$I_{(k-1)}$&$A_{(k)}$\\
\end{tabular}}.
$$

We fix $k$ and $t'$ and use now an induction on $t_1$ such that $t_1\leq k-1$. First, if $t_1=0$, the matrix reduces to:
$$
A_{nk(0,t')}=A_{(k)}=A_{n'kt'},
$$
which is injective.\\

We now suppose that $0<t_1\leq k-1$ and, if $s_1=t_1-1$, $s'=t'$, $A_{nks}$ is injective.

Consider then a principal submatrix in $A_{nkt}$, namely:
$$
P={\tiny\begin{tabular}{c|ccccc}
\diagbox{$a'$}{$a$}&$t_1$&$t_1-1$&$\cdots$&$1$\\
\hline
$t_1$&$A_{(k-t_1)}$&0&$\cdots$&$0$\\
$t_1-1$&$I_{(k-t_1)}$&$A_{(k-t_1+1)}$&&\\
$\vdots$&&$\ddots$&$\ddots$&&\\
%$\vdots$&&&$\ddots$&$\ddots$&&&\\
%$a$&&&&$B_{d_{k-a-1},d_{k-a}}$&$0$&&&\\
%$a-1$&&&&$I_{k-a}$&$B_{k-a+1}$&&&\\
$\vdots$&&&$\ddots$&$\ddots$&\\
%$\vdots$&&&&&&$\ddots$&$\ddots$&\\
$1$&$0$&$\cdots$&$I_{(k-2)}$&$A_{(k-1)}$\\
%$k-|t'|-1$&$0$&$0$&$\cdots$&$0$&$I_{|t'|}$\\
\end{tabular}}.
$$
So that $A_{nkt}=\left[\begin{smallmatrix}P&0\\ R&A_{(k)}\end{smallmatrix}\right]$, with the `row matrix' $R=[\begin{matrix}0&\dots&I_{(k-1)}\end{matrix}]$.\\

Put  now $j=k-1$, %$s=t-\varepsilon_1=(t_1-1,t_2,\ldots,t_n)$. Observe then that
since $s'=t'=(t_2,\ldots,t_n)$, and $j-s_1=k-t_1$, therefore $t_1<k$ is equivalent to $s_1<j$, moreover we have $j<j+1=k\leq |t'|=|s'|$ thus $s_1<j<|s'|$. In other words, the decomposition of the matrix $A_{njs}$ is similar to the decomposition of  $A_{nkt}$, and if $b=a-1$, $(A_{njs})_{b,b}=A_{n'(j-b)s'}=A_{n'(k-a)t'}=A_{(k-a)}$, similarly $I_{(j-b)}=I_{(k-a)}$ hence:
$$
A_{njs}={\tiny\begin{tabular}{c|ccccc}
\diagbox{$b'$}{$b$}&$s_1$&$s_1-1$&$\cdots$&$0$\\
\hline
$s_1$&$A_{(k-t_1)}$&0&$\cdots$&$0$\\
$s_1-1$&$I_{(k-t_1)}$&$A_{(k-t_1+1)}$&&\\
$\vdots$&&$\ddots$&$\ddots$&&\\
%$\vdots$&&&$\ddots$&$\ddots$&&&\\
%$a$&&&&$B_{d_{k-a-1},d_{k-a}}$&$0$&&&\\
%$a-1$&&&&$I_{k-a}$&$B_{k-a+1}$&&&\\
$\vdots$&&&$\ddots$&$\ddots$&\\
%$\vdots$&&&&&&$\ddots$&$\ddots$&\\
$0$&$0$&$\cdots$&$I_{(k-2)}$&$A_{(k-1)}$\\
%$k-|t'|-1$&$0$&$0$&$\cdots$&$0$&$I_{|t'|}$\\
\end{tabular}}=P.
$$

Now we have $j-s_1=k-t_1<\frac{|t'|+1}{2}=\frac{|s'|+1}{2}$, thus we are in the case 2 for $A_{njs}$ and since $s_1=t_1-1$, %is even and $s_2$ is odd, then 
$P=A_{njs}$ is injective. Since, thanks to our induction hypothesis on $n$, $A_{(k)}$ is also injective, this implies that $A_{nkt}$ is injective too.\\

As in the preceding case, we have:
$$
h^2_{nkt}=\dim(\coker(f_{nkt}))=d_{n'(k'-t_1)t'}-d_{n'k't'}.
$$
\end{itemize}

This ends the proof of the theorem.\\ 

\end{proof}

\begin{cor}\label{cor1}
Put $k=\mu-|\lambda|$, $t_i=-2\lambda_i$ and:
$$
d_{nkt}=\#\{\alpha=(\alpha_1,\ldots,\alpha_n):~\alpha_i\leq t_i~(i=1,\dots,n),~|\alpha|=\sum\alpha_i=k\}.
$$
Then for any $n,k,t$, the dimension $h^j_{nkt}$ of the cohomolgy groups $\HD{j}$ ($j=0,1$ or $2$) is given by the following table
$$
%\tiny
{\begin{tabular}{c|cc|ccc|cc}
\small{$k$}&\small{$0$}&\hfill&\tiny{$\frac{|t|+1}{2}$}&\hfill&\tiny{$|t|+1$}&\hfill&\tiny{$+\infty$}\\
\hline
&&&&&&\\
\small{$h^2_{nkt}$}&&\small{$0$}&&\small{$d_{n(k-1)t}-d_{nkt}$}&&\small{$0$}&\\
&&&&&&\\
\small{$h^0_{nkt}$}&&$\left(\begin{smallmatrix}n+k-2\\ k\end{smallmatrix}\right)$&&$\left(\begin{smallmatrix}n+k-2\\ k\end{smallmatrix}\right)$\small{$+d_{n(k-1)t}-d_{nkt}$}%\small{$h^2_{nkt}$}
&&$\left(\begin{smallmatrix}n+k-2\\ k\end{smallmatrix}\right)$&\\
&&&&&&\\
\small{$h^1_{nkt}$}&&$\left(\begin{smallmatrix}n+k-2\\ k\end{smallmatrix}\right)$&&$\left(\begin{smallmatrix}n+k-2\\ k\end{smallmatrix}\right)$\small{$+2(d_{n(k-1)t}-d_{nkt})$}%\small{$h^2_{nkt}$}
&&$\left(\begin{smallmatrix}n+k-2\\ k\end{smallmatrix}\right)$&\\
&&&&&&\\
\end{tabular}.}
$$
\end{cor}

For instance, if $n=1$, we immediately get $\left(\begin{smallmatrix}n+k-2\\ k\end{smallmatrix}\right)=0$, and
$$
d_{1kt}=\left\{\aligned &0&\text{ if }&k\neq t\\
&1&\text{ if }&k=t.\endaligned\right.
$$
Then our table becomes:
$$
\begin{tabular}{c|cc|c|cc}
\small{$k$}&\tiny{$0$}&\hfill&\tiny{$|t|+1$}&\hfill&\tiny{$+\infty$}\\
\hline
&&&&&\\
\small{$h^2_{1kt}$}&&\small{$0$}&\small{$1$}&\small{$0$}&\\
&&&&&\\
\small{$h^0_{1kt}$}&&\small{$0$}&\small{$1$}&\small{$0$}&\\
&&&&&\\
\small{$h^1_{1kt}$}&&\small{$0$}&\small{$2$}&\small{$0$}&\\
&&&&&\\
\end{tabular}.
$$

More generally, we have the following results:
\begin{thm} Suppose that $t%=(t_1,\,t_2,\,\dots,\,t_n)\in\{0,\,1,\,\dots,\,k-1\}
\in\N^n$, and $t_1\leq\ldots\leq t_n<k$. Then 
\begin{itemize}
\item[1.] If $t_1=0$, then $\dim\HD{2}=h^2_{nkt}=h^2_{n'k't'}$. 

We now suppose $t_1>0$,  and put $r=\sup\{i:~t_i=1\}$ ($r=0$ if ($t_1>1$).\\

\item[2.] If $k>|t|+1$, then $h^2_{nkt}=\dim\HD{2}=0$.

If $k=|t|+1$, then $h^2_{nkt}=1$. More precisely, the space $\HD{2}$ is the set of multiples of the 2-cocycle $c$ defined by
\begin{equation}\label{42}
c(X_{h_1},X_{h_2})= \left|\begin{matrix}h'_1&h'_2\\ h''_1&h''_2\end{matrix}\right|\Omega^t.%{(t_1,\dots,t_n)}.
\end{equation} 

\item[3.]  If $k=|t|$, then $h^2_{nkt}=n-1$, a basis of $\HD{2}$ is the set of cocycles $c_i$ ($i<n$), where:
%, for any fixed $i_0$ such that $t_{i_0}>0$, the space $\mathrm{H}^2(\mathfrak{sl}(2),\mathcal{D}_{\underline{\lambda},\mu})$ can be spanned by the 2-cocycles $\omega_i$ defined by:
\begin{equation}
c_i(X_{h_1},X_{h_2})= \left|\begin{matrix}h'_1&h'_2\\ h''_1&h''_2\end{matrix}\right|\Omega^{t-\varepsilon_i}.%{(t_1,\dots,t_n)-\varepsilon_i}.
\end{equation} 

\item[4.]  If $k=|t|-1$, then $h^2_{nkt}=\frac{n(n-1)}{2}-r$,  a basis of $\HD{2}$ is the set of cocycles $c_{ij}$ ($1\leq i<j< n$), and $c_{2i}$ ($r<i< n$) where:%the space $\mathrm{H}^2(\mathfrak{sl}(2),\mathcal{D}_{\underline{\lambda},\mu})$ can be spanned by some 2-cocycles $\omega_{ij}$ of the form:
\begin{equation}\aligned
 c_{ij}(X_{h_1},X_{h_2})&= \left|\begin{matrix}h'_1&h'_2\\ h''_1&h''_2\end{matrix}\right|\Omega^{t-\varepsilon_i-\varepsilon_j}~\hskip2cm~& (i<j<n),\\
c_{2i}(X_{h_1},X_{h_2})&= \left|\begin{matrix}h'_1&h'_2\\ h''_1&h''_2\end{matrix}\right|\Omega^{t-2\varepsilon_i}~\hskip1cm~& (r<i<n),
\endaligned
 \end{equation} 
%$$
%h^2_{nkt}=\dim \HD{2}=\left\{\begin{array}{ll}0
%\quad &\text{if}\quad k>|t|+1\\[4pt]
%1\quad &\text{if}\quad k=|t|+1\\[4pt]
%n-1\quad &\text{if}\quad k=|t|\\[4pt]
%{n(n-1)\over2}-r&\text{if}\quad k=|t|-1,
%\end{array}\right.
%$$
%if $t=(1,\ldots,1,t_{r+1},\ldots,t_n)$ and $t_{r+1}\geq 2$.
%where $r$ the number of $t_i$ and $r$ is the one of $t_i=1$.

%\item[3.] %Moreover, i
%If $k=|t|+1$, then the space $\HD{2}$ is the set of multiples of the 2-cocycle $c$ defined by
 %\begin{equation}\label{42}
%c(X_{h_1},X_{h_2})= \left|\begin{matrix}h'_1&h'_2\\ h''_1&h''_2\end{matrix}\right|\Omega^{(t_1,\dots,t_n)}.
 %\end{equation} 

%\item[4.] If $k=|t|$, then a basis of $\HD{2}$ is the set of cocycles $c_i$ ($i<n$), where:
%, for any fixed $i_0$ such that $t_{i_0}>0$, the space $\mathrm{H}^2(\mathfrak{sl}(2),\mathcal{D}_{\underline{\lambda},\mu})$ can be spanned by the 2-cocycles $\omega_i$ defined by:
%\begin{equation}
%c_i(X_{h_1},X_{h_2})= \left|\begin{matrix}h'_1&h'_2\\ h''_1&h''_2\end{matrix}\right|\Omega^{(t_1,\dots,t_n)-\varepsilon_i}.
%\end{equation} 
 
%\item[5.] If $k=|t|-1$, then  a basis of $\HD{2}$ is the set of cocycles $c_{ij}$ ($1\leq i<j\leq n$), and $c_{2i}$ ($r<i\leq n$) where:%the space $\mathrm{H}^2(\mathfrak{sl}(2),\mathcal{D}_{\underline{\lambda},\mu})$ can be spanned by some 2-cocycles $\omega_{ij}$ of the form:
%\begin{equation}
 %\omega_{ij}(X_{h_1},X_{h_2})= \left|\begin{matrix}h'_1&h'_2\\ h''_1&h''_2\end{matrix}\right|\Omega^{(t_1,\dots,t_n)-\varepsilon_i-\varepsilon_j}.
 %\end{equation} 
\end{itemize}

\end{thm}
\begin{proof}
\
\begin{itemize}
\item[1.] If $k>|t|+1$ then $K^{k-1}_t=0$. 

If $k=|t|+1$ then $K^{k-1}_t$ is 1-dimensional spanned by $\Omega^t%{(t_1,\dots,t_n)}
$ and $ K^{k}_t=0$. In fact, the proof of Lemma \ref{fcohomology} shows that \eqref{42} is a basis of the space $\HD{2}$.\\%\mathrm{H}^2(\mathfrak{sl}(2),\mathcal{D}_{\underline{\lambda},\mu})$.
   
\item[2.] If $k=|t|$ then $K^{k-1}_t$ is $n$-dimensional, with basis $\{E_{t-\varepsilon_i}\}$ (see \eqref{Eformula} for the definition of $E_\a$) and  $K^{k}_t$ is $n$-dimensional. But 
$$
f(E_{t})=\sum_{i=1}^nE_{t-\varepsilon_i},
$$
we can therefore identify $\coker(f)$ to the subspace in $K^{k-1}_t$%, supplementary to $f(K^k_t)$, spanned by the 
with basis $\{\Omega^{%(t_1,\dots,t_n)
t-\varepsilon_i}:~i<n\}$. We conclude by using the proof of Lemma \ref{fcohomology}.\\

\item[3.] If $k=|t|-1$ then $K^{k-1}_t$ is $({n(n-1)\over2}+n-r)$-dimensional, with  basis:% generated by 
 $$
 \big\{E_{t-\varepsilon_i-\varepsilon_j},~~1\leq i<j\leq n\big\}\cup\big\{E_{t-2\varepsilon_i},~~r< i\leq n\big\}
 $$
 and  $K^{k}_t$ is $n$-dimensional with basis $\{E_{t-\varepsilon_i},~~1\leq i\leq n\}$.
 
The space $f(K^{k}_t)$ is spanned by the
$$
f(E_{t-\varepsilon_i})=\sum_{j\neq i}E_{t-\varepsilon_i-\varepsilon_j}+E_{t-2\varepsilon_i}
$$
(the second term happens only if $i>r$). Using the inverse of the lexigographic ordering, the vector with the smallest index in the decomposition of $f(E_{t-\varepsilon_i})$ is $E_{t-\varepsilon_i-\varepsilon_n}$ ($i\leq n$), hence the vectors $f(E_{t-\varepsilon_i})$ are independent, and $h^2_{nkt}=\dim\coker(f)$ is $\frac{n(n-1)}{2}-r$.

Moreover, the vectors $\O^{t-\varepsilon_i-\varepsilon_j}$ ($1\leq i<j<n$) and $\O^{t-2\varepsilon_i}$ ($r<i<n$) form a basis for a subspace, supplementary to $f(K^k_t)$ in $K^{k-1}_t$. With Lemma \ref{cohomology}, this achieves the proof of the theorem.\\
\end{itemize}
\end{proof}

%-----------------------------------------------------------------------------------------------------------------------

\section{Small $n$ cases}

\

%-----------------------------------------------------------------------------------------------------------------------

To illustrate our result, let us consider the small $n$ cases, computing the $h^j_{nkt}$ for $n=1,2,3,4$.
The case $n=1$ was presented in the preceding section. 
Let us consider now the cases $n>1$. We still suppose that $t_1\leq  \dots \leq t_n$.  Recall that $d_{nkt}=\# D_{nkt}$ and if $n'=n-1$ and $t'=(t_2,\ldots,t_n)$, hence, putting $j=k-a$, and since, by definition, $d_{nkt}=0$ if $k<0$ or $k>|t|$, 
\begin{equation}\label{rec}\aligned
d_{nkt}&=\sum_{a=\sup(0,k-|t'[)}^{\inf(t_1,k)}d_{n'(k-a)t'}=\sum_{j=\sup(0,k-t_1)}^{\inf(k,|t'|)}d_{n'jt'}.
%&=\sum_{j=k-t_1}^kd_{n'jt'}.
\endaligned
\end{equation}
Let us define $x_{nkt}$ by:
$$
x_{nkt}=d_{nkt}-d_{n(k-1)t}.
$$
(Hence $x_{nkt}=0$ if $k<0$ or $k>|t|+1$, $x_{n0t}=1$.) It is clear that we also have:
$$
d_{nkt}=\sum_{j=0}^k x_{njt},\text{ and }h^2_{nkt}=\sup(0,-x_{nkt}).
$$

Therefore, we have induction property for $x_{nkt}$:
\begin{equation}\label{induction}
\aligned
x_{nkt}&=\sum_{j=k-t_1}^{k}d_{n'jt'}-\sum_{j=k-t_1-1}^{k-1}d_{n'jt'}\\
&=%\left\{\aligned 
%&d_{n'kt'}=\sum_{j=0}^kx_{n'jt'},&\text{ if } k\leq t_1,\\
d_{n'kt'}-d_{n'(k-1-t_1)t'}=\sum_{j=k-t_1}^kx_{n'jt'}.%,&\text{ if } t_1< k\leq |t'|,\\
%&-d_{n'(k-1-t_1)t'}=-\sum_{j=0}^{k-1-t_1}x_{n'jt'},&\text{ if } |t'|< k\leq |t|+1.
% \endaligned\right.
\endaligned
\end{equation}

On the other hand,  the matrix $A_{n(|t|+1-k)t}$ has same size that the matrix $A_{nkt}^t$, namely is a $d_{nkt}\times d_{n(k-1)t}$-matrix, then
$$
d_{n(|t|+1-k)t}=d_{n(k-1)t},~~x_{n(|t|+1-k)t}=d_{n(|t|+1-k)t}-d_{n(|t|-k)t}=-x_{nkt}.
$$
Therefore, we can restrict ourselves to the computation of the $x_{nkt}$ for $k\leq \left[\frac{|t|+1}{2}\right]$. Since $x_{1kt}=0$ except for $k=0$ or $|t|+1$, where $x_{10t}=1=-x_{1(|t|+1)t}$, we get for $n=2$ or $3$:

\subsection{Case n=2.} From the preceding relations, we immediatly get the $x_{2kt}$ and the $h^2_{2kt}$. These results are well known, and we leave the computation to the reader.\\

\begin{lem}\label{n=2}

The $x_{2kt}$ and the dimension $h^2_{2kt}$ of the group $H^2(\mathfrak{sl}(2),\Dlm)$ are given by the table:
$$
\begin{tabular}{c|ccc|ccc|ccc|ccc}
\small{$k$}&\small{$0$}&\hfill&\small{$t_1$}&\small{$t_1+1$}&\hfill&\small{$t_2$}&\small{$t_2+1$}&\hfill&\small{$|t|+1$}&\small{$|t|+2$}\hfill&\small{$+\infty$}\\
\hline
\small{$x_{2kt}$}&&\small{$1$}&&&\small{$0$}&&&\small{$-1$}&&&\small{$0$}&\\
\hline
\small{$h^2_{2kt}$}&&\small{$0$}&&&\small{$0$}&&&\small{$1$}&&&\small{$0$}&\\
\end{tabular}
$$
\end{lem}

\subsection{Case n=3.} Similarly, we can compute the $x^2_{3kt}$ and the $h^2_{3kt}$, using formula \eqref{induction}. But here the result depends of the relative position of $t_1+t_2$ and $t_3$. We get:\\

\begin{lem}\label{n=3}

If $n=3$, the dimensions $h^2_{2kt}$ of the group $\HD{2}%H^2(\mathfrak{sl}(2),\Dlm)
$ are given by the table, where $a=\left[\frac{|t|+1}{2}\right]+1$, $t_{ij}=t_i+t_j$ and $b=\sup(t_{12},t_3)$:
$$
\tiny{\begin{tabular}{c|c|c|c|c}
$k$&$a$\hskip 3cm$b$&$b+1$\hskip 1cm$t_{13}$&$t_{13}+1$\hskip 1cm$t_{23}$&$t_{23}+1$\hskip 1cm$|t|+1$\\
\hline
$h^2_{3kt}$&$\left\{\aligned &0&&(t_{12}<t_3)\\ &2k-1-|t|&&(t_3<t_{12}) \endaligned\right.$&$k-t_3$&$t_1+1$&$|t|+2-k$\\
\end{tabular}}
$$

\end{lem}

Some of these results were known (see \cite{be}).\\

\subsection{Case n=4.} Looking to the $x_{4kt}$ for $k\leq\left[\frac{|t|+1}{2}\right]$, We see immediatly that, if $k\leq t_1$,
$$
x_{4kt}=\sum_{j=0}^k x_{3j(t_2t_3t_4)}=\sum_{j=0}^k(j+1)=\frac{1}{2}(k+2)(k+1),
$$
and if $k\geq t_2+t_3+t_4+1=|t|+1-t_1$:
$$
x_{4kt}=-x_{4(|t|+1-k)t}=-\frac{1}{2}(|t|-k+3)(|t|-k+2).
$$
Suppose now that $k>t_1$, then $t_1+1\leq k\leq \frac{t_1+1+|t'|}{2}$, therefore $k\leq|t'|=t_2+t_3+t_4$, and:
$$
x_{4kt}=\sum_{j=k-t_1}^k x_{3j(t_2t_3t_4)}.
$$
To determine this sum, we have to compare $k-t_1$ and $k$ with $0$, $t_2$, $t_3$, $t_2+t_3$, $t_4$, $t_2+t_4$, $t_3+t_4$ and $t_2+t_3+t_4$. Now let us write $t_{ab}$ for $t_a+t_b$, $t_{abc}$ for $t_a+t_b+t_c$, we have the following partial ordering for these numbers:
$$
\begin{array}{ccccc}
&&t_{23}&&\\
&\swarrow&&\nwarrow&\\
0\leftarrow t_2 \leftarrow t_3 &&&& t_{24}\leftarrow t_{34}\leftarrow t_{234}.\\
&\nwarrow&&\swarrow&\\
&&t_{4}&&\\
\end{array}
$$
Suppose we insert $k-t_1$ in this partially ordered set, we then have to insert $k$ in the set: 
$$
\begin{array}{ccccc}
&&t_{123}&&\\
&\swarrow&&\nwarrow&\\
t_1\leftarrow t_{12} \leftarrow t_{13} &&&& t_{124}\leftarrow t_{134}\leftarrow t_{1234}.\\
&\nwarrow&&\swarrow&\\
&&t_{14}&&\\
\end{array}
$$

Finally, we have to insert $k$ in the partially ordered set:
$$\tiny{
\begin{array}{ccccccccccccc}
&&t_{12}&\leftarrow&t_{13}&\leftarrow&t_{23}&\leftarrow&t_{123}&\leftarrow&t_{124}&&\\
&\swarrow&&&&&&&&&&\nwarrow&\\
0\leftarrow t_1 \leftarrow t_2 &&&\swarrow&&\nwarrow&&\nwarrow&&\swarrow&&& t_{134}\leftarrow t_{234}\leftarrow t_{1234}.\\
&\nwarrow&&&&&&&&&&\swarrow&\\
&&t_{3}&\leftarrow&t_{4}&\leftarrow&t_{14}&\leftarrow&t_{24}&\leftarrow&t_{34}&&\\
\end{array}
}
$$
It is possible to list all the 10 possible total orderings extending this partial ordering. For each of these orderings, we can then present a table ginving the $x_{4kt}$ for $k$ in each of the 13 places between $t_2$ and $t_{234}$. For instance, suppose the total ordering is:
$$\tiny{\begin{array}{c}
0\leftarrow t_1\leftarrow t_2\leftarrow t_{12}\leftarrow t_3\leftarrow t_{13}\leftarrow t_4\leftarrow t_{23}\leftarrow t_{14}\leftarrow t_{123}\leftarrow t_{24}\leftarrow t_{124}\leftarrow t_{34}\leftarrow t_{134}\leftarrow t_{234}\leftarrow t_{1234},
\end{array}}
$$
then the table is:
$$\tiny{
\begin{tabular}{c|c}
$k$&$x_{4kt}$\\
\hline
$0$&\\
&$\frac{1}{2}k^2+\frac{3}{2}k+1$\\
$t_1$&\\
&$(1+t_1)k+(1+t_1)(1-\frac{1}{2}t_1)$\\%(-\frac{1}{2}t_1+k+1)(1+t_1)$\\
$t_2$&\\
&$-\frac{k^2}{2}+(t_1+t_2+\frac{1}{2})k-\frac{1}{2}t_1(t_1-1)-\frac{1}{2}t_2(t_2-1)+1$\\
$t_{12}$&\\
&$(t_1+1)(t_2+1)$\\
$t_3$&\\
&$-\frac{k^2}{2}+\frac{2t_3-1}{2}k-\frac{t_3}{2}(t_3+2)+(t_1+1)(t_2+1)$\\%(1/2)*t3^2+(1/2)*(-2*t3+1)*t3-(1/2)*k^2+(1/2)*(2*t3-1)*k+(t2+1)*(1+t1)$\\
$t_{13}$&\\
&$-(t_1+1)(k-\frac{t_1}{2}-t_2-t_3-1)$\\
$t_4$&\\
&$-\frac{k^2}{2}-(t_1-t_4+\frac{3}{2})k+\frac{t_1^2-t_4^2}{2}+(t_1+1)(t_2+t_3)+\frac{t_4}{2}+\frac{3t_1}{2}+1$\\
$t_{23}$&\\
&$\frac{1}{2}(t_{1234}+1-2k)(\frac{t_{123}-t_4}{2}+1)$\\
$t_{14}$&\\
&$\frac{k^2}{2}-(2t_{123}+\frac{5}{2})k+(t_{1234}+2)t_1+(\frac{t_2}{2}+t_3+\frac{1}{2})t_2+(\frac{t_3}{2}+\frac{3}{2})t_3+t_4+1$\\
$t_{123}$&\\
%\end{tabular}}
%$$
%$$\tiny{
%\begin{tabular}{c|c}
%$k$&$x_{4kt}$\\
%\hline
%$t_{123}$&\\
&$(t_1+1)(\frac{t_1}{2}+t_4-k) $\\
$t_{24}$&\\
&$\frac{k^2}{2}+(t_{124}+\frac{3}{2})(t_4-k)+\frac{t_1}{2}(t_1+1)+\frac{t_2}{2}(t_2+1)$\\
$t_{124}$&\\
&$-(t_1+1)(t_2+1)$\\
$t_{34}$&\\
&$\frac{k^2}{2}-(t_ {34}+\frac{1}{2})k+(t_3+\frac{t_4}{2}+\frac{1}{2})t_4+(t_3+\frac{1}{2})t_3-\frac{t^2}{2}-(t_2+1)(t_1+1)$\\
$t_{134}$&\\
&$(t_1+1)(k-t_{1234}+\frac{t_1}{2}-2)$\\
$t_{234}$&\\
&$-\frac{1}{2}(t_{1234}+2-k)(t_{1234}+3-k)$\\
$t_{1234}+1$&\\
\end{tabular}}
$$

%-----------------------------------------------------------------------------------------------------------------------

\section{Homogeneous case}\label{homogeneous}

\

%-----------------------------------------------------------------------------------------------------------------------

In this section, we restrict ourselves to the case where the $t_i$ are equal. Let us call this case the homogeneous case. In this section, we put $t=t_i$ for $i=1,\dots,n$ and write $d_{nkt}$, $x_{nkt}$ and so on for $d_{nk(t,t,\ldots,t)}$, $x_{nk(t,t,\ldots,t)}$. In this very simple case, it is possible to study the maps $k\mapsto x_{nkt}$ and $k\mapsto h^2_{nkt}$, for $n>1$.\\

%With the notation of the theorem, if $r=n$, it is also easy to compute the $h^2_{nkt}$ for any $k\leq n$:\\

First, let us consider the simplest case, when $t=1$.\\
 
\begin{prop} 
Suppose $t=1$. Then
%If $t_1=\ldots=t_n=1$ and $k=n-m$, then
$$
%\mathrm{dim}\mathrm{H}^2_\mathrm{diff}(\mathfrak{sl}(2),\mathcal{D}_{\underline{\lambda},\mu})
h^2_{nk1}=\dim\HD{2}=\sup\left(\frac{(%n-2m-1
2k-1-n)n!}{%(n-m)
k!(%m
n-k+1)!},0\right).
$$ 
Moreover there is an explicit recursive construction of a set of indexes $\mathcal A_{nk}$ such that a basis of $\HD{2}$ is given by the cocycles:
$$
c_\a(X_{h_1},X_{h_2})= \left|\begin{matrix}h'_1&h'_2\\ h''_1&h''_2\end{matrix}\right|\Omega^\a\hskip 1,5cm (\a\in\mathcal A_{nk}).
$$
\end{prop}

\begin{proof}
In fact, with our hypothesis, for any $s$,  if $\a_X$ is the characteristic function of the subset $X$, the set
$$
\{E_{\a_X}:~X\subset\{1,\ldots,n\},~\#X=s\}
$$
is a basis of $K^s_\1$% is $\{E_{\a_X}:~X\subset\{1,\ldots,n\},~\#X=s\}$, if $\a_X$ is the characteristic function of the subset $X$
. Hence $\dim K_1^s=\left(\begin{smallmatrix}n\\ s\end{smallmatrix}\right)$. Therefore:
$$
\dim(K_1^{k-1})-\dim(K^k_1)=\left(\begin{matrix}n\\ k-1\end{matrix}\right)-\left(\begin{matrix}n\\ k\end{matrix}\right)=\frac{(2k-1-n)n!}{k!(n-k+1)!}.
$$

Observe that if $k=n$, the matrix $A_{nk1}$ is a column matrix with $k$ rows, each being 1. Therefore a basis of $\coker(A_{kk1})$ is the set of $\O^\a$, with $\alpha_i=1$ except for an unique $i<k$, its dimension is $k-1$.

If $n=2k-1$, the matrix $A_{nk1}$ is invertible, $\coker(A_{(2k-1)k1})=0$.

Now if $k<n<2k-1$ then $A_{nk1}$ has the form:
$$
A_{nk1}=\left[\begin{matrix}A_{(n-1)(k-1)1}&0\\ I&A_{(n-1)k1}\end{matrix}\right].
$$
And it is easy to prove that:
$$
\coker(A_{nk1})=\coker(A_{(n-1)(k-1)1})\oplus\coker(A_{(n-1)k1}).
$$

The entries in $A_{(n-1)(k-1)1}$ are the $x_{\a,\a'}$ such that $\a_1=\a'_1=1$, and the entries in $A_{n-1)k1}$ are the $x_{\a,\a'}$ such that $\a_1=\a'_1=0$. We write in both cases $\a=\beta\gamma$ with $\beta=(1)$ (resp. $(0)$), $\gamma=(\a_2,\ldots,\a_n)$, and $\a'=\beta\gamma'$.

If $n-1<2(k-1)-1$, we repeat this decomposition for $A_{(n-1)(k-1)1}$, if $k<n-1$, we decompose $A_{(n-1)k1}$. The process stops when, for each matrix $A_{mj1}$ in the decomposition, either $m=2j-1$ or $m=j$. Since $\coker(A_{nk1})=\sum\coker(A_{mj1})$, this gives a basis of $\coker(A_{nk1})$.

\end{proof}

Let us consider the example $n=7$, $k=5$. Let $p$ be the number of  the iteration step, in the following table, we put one star if $(m,j)=(2j-1,j)$ and two stars if $(m,j)=(j,j)$.
$$
\begin{array}{c|c|c|c|c|c}
%\hline
p&0&1&2&3&4\\
\hline
\beta&(~)&\begin{array}{c}(1)\\ (0)\end{array}&\begin{array}{c}(11)^*\\ (10)\\ (01)\\ (00)^{**}\end{array}&\begin{array}{c}(101)\\ (100)^{**}\\ (011)\\ (010)^{**}\end{array}&\begin{array}{c}(1011)^*\\ (1010)^{**}\\ (0111)^*\\ (0110)^{**}\end{array}\\
\hline
(m,j)&(7,5)&\begin{array}{c} (6,4)\\ (6,5)\end{array}&\begin{array}{c} (5,3)^*\\ (5,4)\\ (5,4)\\ (5,5)^{**}\end{array}&\begin{array}{c} (4,3)\\ (4,4)^{**}\\ (4,3)\\ (4,4)^{**}\end{array}&\begin{array}{c} (3,2)^*\\ (3,3)^{**}\\ (3,2)^{*}\\ (3,3)^{**}\end{array}\\
\hline
\tiny{\a= \beta|\gamma}&&&\tiny{\begin{array}{c} 00|11101\\ 00|11011\\ 00|10111\\ 00|01111\end{array}}&\tiny{\begin{array}{cc} 100|1101&010|1101\\ 100|1011& 010|1011\\ 100|0111& 010|0111\end{array}}&\tiny{\begin{array}{cc} 1010|101&0110|101\\ 1010|011&0110|011\end{array}}\\
%\hline
\end{array}
$$
Let $\mathcal A_{75}$ be the set of the 14 $\a$ in the last row of this table. Then a basis of $\coker(A_{751})$ is $\{\O^\a:~\a\in\mathcal A_{75}\}$.\\

Coming back to the general homogeneous case, we first recall that:
$$
d_{nkt}=\#\{\alpha\in\N^n:~\alpha_i\leq t,~|\alpha|=k\},\hskip 1cm x_{nkt}=d_{nkt}-d_{n(k-1)t},
$$ 
and if $n'=n-1$, $t'=(t,\ldots,t)\in \N^{n'}$,
$$
x_{nkt}=\sum_{j=\sup(0,k-t)}^{\inf(k,n't)}x_{n'jt'}, \hskip 1cm h^2_{nkt}=\sup(-x_{nkt},0).
$$
%Put also
%$$
%b_{n,k}=\#\{\alpha\in\N^n:~|\alpha|=k\}=\left(\begin{matrix} n+k-1\\ k\end{matrix}\right).
%$$

We saw that $x_{nkt}=0$ if $k>nt+1$. Let us suppose now $n>1$ and decompose the segment $[0,nt+1]$ into $n$ sub-segments, by putting:
$$\aligned
I_0&=[0,t],\\
I_p&=[pt+1,(p+1)t]\hskip 0.5cm (p=1,2,\dots ,n-2),\\
I_{n-1}&=[(n-1)t+1,nt+1].
\endaligned
$$

Recall that in the preceding section, we saw that:
$$
\begin{tabular}{c|c|c}
\small{$I_p$}&\small{$I_0$}&\small{$I_1$}\\
\hline
\small{$x_{2kt}$}&\small{$1$}&\small{$-1$}\\
\end{tabular},\hskip 0.7cm
\begin{tabular}{c|c|c|c}
\small{$I_p$}&\small{$I_0$}&\small{$I_1$}&\small{$I_2$}\\
\hline
\small{$x_{3kt}$}&\small{$k+1$}&\small{$-2k+2t+1$}&\small{$k-3t-2$}\\
\end{tabular}.
$$
Moreover, if $n=4$, a direct computation gives that on each $I_p$, the function  $k\mapsto x_{4kt}$ is a quadratic function, namely:
$$\aligned
&\text{on $I_0$, }~~x_{4kt}=\frac{1}{2}k^2+\frac{3}{2}k+1,\\
&\text{on $I_1$, }~~x_{4kt}=-\frac{3}{2}k^2+\frac{1}{2}(3t-1)k+(t^2+2t-1),\\
&\text{on $I_2$, }~~x_{4kt}=\frac{3}{2}k^2-7(t+\frac{1}{2})k+(7t^2+7t+1),\\
&\text{on $I_3$, }~~x_{4kt}=-\frac{1}{2}k^2+\frac{1}{2}(8t+5)k-(8t^2+10t+\frac{5}{2}).\\
\endaligned
$$

Let us now prove:\\

\begin{prop}\label{homogene case}

With our notation, on each $I_p$, the function $k\mapsto x_{nkt}$ is a polynomial function, with degree $n-2$, and the term in $k^{n-2}$ is:
$$
\frac{(-1)^p}{(n-2)!}\left(\begin{matrix}n-1\\ p\end{matrix}\right)~k^{n-2}
$$
\end{prop}

\begin{proof}
Recall the Faulhaber formula if $B_r$ is the $r^{th}$ Bernouilli number,
$$
\sum_{j=0}^kj^m=\frac{1}{m+1}\big(k^{m+1}+\frac{1}{2}(m+1)k^m+\sum_{r=2}^m \left(\begin{matrix}m+1\\ r\end{matrix}\right)B_r k^{m+1-r} \big).
$$

The proposition holds for $n=2,~3$ or 4. Suppose it is holding for some $n'=n-1\geq 2$, then if $k\leq t$ (or $p=0$) (resp. $k>n't$ or $p=n-1$),
$$\aligned
x_{nkt}&=\sum_{j=0}^k x_{n'jt}=\sum_{j=0}^k\big(\frac{1}{(n'-2)!}j^{n'-2}+\sum_{s<n'-2}a_s(t)j^s\big)\\
&=\frac{1}{(n-2)!}k^{n-2}+\sum_{s<n-2}b_s(t)k^s.
\endaligned
$$
Respectively:
$$\aligned
x_{nkt}&=-x_{n(nt+1-k)t}=-\frac{1}{(n-2)!}(nt+1-k)^{n-2}+\sum_{s<n-2}a_s(t)k^s\\
&=\frac{(-1)^{n-1}}{(n-2)!}k^{n-2}+\sum_{s<n-2}b_s(t)k^s.
\endaligned
$$
Suppose now $k\in I_p$, with $0<p<n-1$, then:
$$\aligned
x_{nkt}&=\sum_{j=k-t}^{pt}x_{n'jt}+\sum_{j=pt+1}^k x_{n'jt}\\
&=\sum_{j=k-t}^{pt}\frac{(-1)^{p-1}}{(n'-2)!}\left(\begin{matrix}n'-1\\ p-1\end{matrix}\right)j^{n'-2}+\sum_{j=pt+1}^k\frac{(-1)^{p}}{(n'-2)!}\left(\begin{matrix}n'-1\\ p\end{matrix}\right)j^{n'-2}+\\
&\hskip 3cm+\sum_{s<n'-2}a_s(t)j^s\\
%&=-\sum_{j=0}^{k-t-1}\frac{(-1)^{p-1}}{(n'-2)!}\left(\begin{matrix}n'-1\\ p-1\end{matrix}\right)j^{n'-2}+\sum_{j=pt+1}^k\frac{(-1)^{p}}{(n'-2)!}\left(\begin{matrix}n'-1\\ p\end{matrix}\right)j^{n'-2}+\\
%&\hskip 3cm+\sum_{j=0}^{pt}\frac{(-1)^{p-1}}{(n'-2)!}\left(\begin{matrix}n'-1\\ p-1\end{matrix}\right)j^{n'-2}+\sum_{s<n'-2}a_s(t)j^s\\
%&=\frac{(-1)^p}{(n-2)!}\big(\left(\begin{matrix}n'-1\\ p-1\end{matrix}\right)+\left(\begin{matrix}n'-1\\ p\end{matrix}\right)\big)k^{n-2}+\sum_{s<n-2}b_s(t)k^s\\
&=\frac{(-1)^p}{(n-2)!}\left(\begin{matrix}n-1\\ p\end{matrix}\right)k^{n-2}+\sum_{s<n-2}b_s(t)k^s.
\endaligned
$$
This proves the proposition.\\

\end{proof}

%-----------------------------------------------------------------------------------------------------------------------

%\section{Bibliography}\label{bibliography}

%-----------------------------------------------------------------------------------------------------------------------

%----------------------------------------------------------------------------------------------------------------------- 

%\input{Introduction}

%\input{Various cohomologies}

%\input{Map F}

%\input{Generic}

%\input{from F to f}

%\input{maximalrank}

%\input{small n cases}

%\input{homogene}

%\input{Bibliography}

%-----------------------------------------------------------------------------------------------------------------------

\end{document}